\documentclass{amsart}
\usepackage[utf8]{inputenc}
\usepackage{amsmath}
\usepackage{amsfonts}
\usepackage{biblatex}
\usepackage{amssymb}

\usepackage{textcomp}
\usepackage{gensymb}
\usepackage{pst-node}
\usepackage{auto-pst-pdf}
\usepackage{tikz-cd}
\usepackage{tikz,tkz-euclide}
\usepackage{siunitx}
\usepackage{pgfplots}
\pgfplotsset{compat=1.17}
\usepackage{float} 
\bibliography{bib}

\usepackage{amsthm}

\theoremstyle{plain}
\newtheorem{prop}{Proposition}[section]

\newtheorem{thm}{Theorem}[section]
\newtheorem{lem}{Lemma}[section]

\theoremstyle{definition}
\newtheorem{eg}{Example}[section]

\newtheorem{defn}{Definition}[section]

\newcommand{\Z}{\mathbb Z}
\newcommand{\F}{\mathbb F}
\newcommand{\ModCmp}[3]{#1 \equiv #2 \;(\bmod\; #3)}

\newcommand{\SL}{\mathrm{SL}_3(\F_7)}
\newcommand{\PSL}{\mathrm{PSL}_3(\F_7)}
\newcommand{\GL}{\mathrm{GL}_3(\F_7)}
\newcommand{\GLLarge}{\mathrm{GL}_3(\F_{7^3})}
\newcommand{\SLTwo}{\mathrm{SL}_2(\F_7)}
\newcommand{\GenSL}{\mathrm{SL}_n(\F_7)}
\newcommand{\GLTwo}{\mathrm{GL}_2(\F_7)}
\newcommand{\GenGL}{\mathrm{GL}_n(\F_7)}

\newcommand{\MatrThree}[9]{
    \begin{pmatrix}
        #1&#2&#3\\
        #4&#5&#6\\
        #7&#8&#9\\
    \end{pmatrix}
}

\newcommand{\vecthree}[3]{
    \begin{pmatrix}
        #1\\
        #2\\
        #3\\
    \end{pmatrix}
}

\newcommand{\vectwo}[2]{
    \begin{pmatrix}
        #1\\
        #2\\
    \end{pmatrix}
}

\newcommand{\MatrTwo}[4]{
    \begin{pmatrix}
        #1&#2\\
        #3&#4\\
    \end{pmatrix}
}

\DeclareMathOperator{\lcm}{lcm}
\DeclareMathOperator{\tr}{tr}

\newcommand{\mbf}[1]{\mathbf{#1}}

\newcommand{\vv}[1]{\boldsymbol{#1}}

\title{Eigenvalue-less matrices from $\SL$ to $\PSL$}
\author[J.L.C,G.C,Y.J, L.L, J.L, W.L, M.S, H.W]{Juan Lucas Callo, George Chen, Yasiru Jayasooriya*,  Leo Li, Jingni Liao, William Liu, Michael Sun, Haibing Wang}
\date{October 2021}

\begin{document}

\begin{abstract}
We analyse the set of matrices in $\SL$ without eigenvalues explicitly, extracting nice bijections between the 18 equally sized conjugacy classes contained within. In doing so, we discover a set of $18$ commuting matrices for which every conjugacy class is represented and tells us how to decide when collections of commuting matrices are simultaneously conjugate. The main innovation is the proofs are accessible to undergraduates and do not rely on computer calculations. It is also rare that the details of special cases are written down.
\end{abstract}

\maketitle
\section*{Introduction}

 The \emph{special linear group} $\SL$ consists of $3 \times 3$ matrices
 $$
 \MatrThree{a}{b}{c}{d}{e}{f}{g}{h}{i}.
 $$
with entries $a,b,c,d,,,,i\in\Z/7\Z$, and determinant 1, under matrix multiplication.

The size of the group is
$$2^5\times3^3\times 7^3\times 19.$$
Consider the subgroup $Z$ of $\SL$ consisting of the following matrices.

$$
\MatrThree{1}{0}{0}{0}{1}{0}{0}{0}{1},
\MatrThree{2}{0}{0}{0}{2}{0}{0}{0}{2},
\MatrThree{4}{0}{0}{0}{4}{0}{0}{0}{4}.
$$

The matrices in this $Z$ commute with all matrices in $\SL$ and therefore form a normal subgroup. The quotient group $\SL / Z$ is known as the projective special linear group or $\PSL$ and is a finite simple group of order
$$2^5\times3^2\times 7^3\times 19.$$

One reason why one may be interested in these groups is that these are small examples for which the Inverse Galois Problem is unknown. Or one might simply be interested to explore how linear algebra works over a finite field with opportunities to also apply one's understanding of the theory of finite groups.

The main goal of our project was to investigate the subset of matrices that do not have any eigenvectors. Some questions that motivated our exploration of these matrices were:

\begin{itemize}
    \item How many matrices do not have eigenvectors in $\SL$?
    \item What orders do these matrices have?
    \item How many conjugacy classes are there?
    \item Can we find nice representatives for each conjugacy class?
    \item How are the conjugacy classes related to each other?
    \item Is it possible to find a set of commuting matrices without eigenvectors where each represents a different conjugacy class?
    \item When are collections of commuting matrices without eigenvectors simultaneously conjugate?
    \item What maximal subgroups cover this set?
\end{itemize}

\begin{thm}\label{ConjClassOverview} 
In $\SL$:

\begin{itemize}
    \item There are exactly $18$ conjugacy classes in $\SL$ where none of the matrices in each conjugacy class have an eigenvector with entries in $\F_7$.
    \item Each of these conjugacy classes has size $2^5 \times 3^2 \times 7^3$.
    \item These $18$ conjugacy classes are the sets $\mbf{[0, 1]}$, $\mbf{[0, 2]}$, $\mbf{[0, 4]}$, $\mbf{[1,0]}$, $\mbf{[1, 3]}$, $\mbf{[1, 6]}$, $ \mbf{[2, 0]}$, $\mbf{[2, 5]}$, $\mbf{[2, 6]}$, $\mbf{[3, 1]}$, $\mbf{[3, 4]}$, $\mbf{[4, 0]}$, $\mbf{[4, 3]}$, $\mbf{[4, 6]}$, $\mbf{[5, 1]}$, $\mbf{[5, 2]}$, $\mbf{[6, 2]}$ and $\mbf{[6, 4]}$.
    \item All the elements in each of the conjugacy classes $\mbf{[0, 2]}$, $\mbf{[1, 3]}$, $\mbf{[2, 0]}$, $\mbf{[3, 1]}$, $\mbf{[3, 4]}$ and $\mbf{[4, 3]}$ have order $19$.
    \item All the elements in each of the conjugacy classes $\mbf{[0, 1]}$, $\mbf{[0, 4]}$, $\mbf{[1, 0]}$, $\mbf{[1, 6]}$, $\mbf{[2, 5]}$, $\mbf{[2, 6]}$, $\mbf{[4, 0]}$, $\mbf{[4, 6]}$, $\mbf{[5, 1]}$, $\mbf{[5, 2]}$, $\mbf{[6, 2]}$ and $\mbf{[6, 4]}$ have order $57$.
\end{itemize}
\end{thm}

Undoubtedly much is known about this group either as a special case of many general theorems or as the result of countless computer calculations whose results are stored in online databases. One of our goals is to look at this single group explicitly on its own and prove all of our results by hand and conceptually to produce a more accessible introduction into the topic for younger students. We have also included an extra section where explanatory notes and omitted details are given to help young students, which takes advantage of the relaxing page restrictions due to electronic publications.

An unexpected outcome of this project was to make observations about collections of commuting matrices and their simultaneous conjugacy.

\begin{thm}
There exists $18$ commuting matrices with each from a different conjugacy class $\mbf{[i, j]}$.
\end{thm}



\begin{thm}
Every collection of commuting matrices excluding $2I,4I$ either consists of matrices all with eigenvectors or they consist of matrices all without eigenvectors and all powers of some $M$.
Two such collections are simultaneously conjugate if and only if there exists $M_1$ for the first collection to be powers of $M_1$ and $M_2$ of the same order for the second collection in which all the corresponding powers agree modulo their common order.
\end{thm}


\section{Irreducible characteristic polynomials over $\F_7$}

We compute the characteristic polynomials of matrices without eigenvectors. These will necessarily be irreducible in $\F_7$.
Let
$$
M := \MatrThree{a}{b}{c}{d}{e}{f}{g}{h}{i} \in \SL
$$
and consider its characteristic polynomial
$$\det(\lambda I-M) = \lambda^3-(a+e+i)\lambda^2+(ae+ie+ai-cg-bd-fh)\lambda+\det(-M).$$
With $\tr(M)=a+e+i$ and $\det(-M) = -\det(M) = -1$, we need the equation
$$\lambda^3-\tr(M)\lambda^2+(ae+ie+ai-cg-bd-fh)\lambda-1=0$$
to not hold for $\lambda=0,1,2,\dots,6$.
That is,
\begin{align*} 
\lambda &= 1: \hspace{20px} -(a+e+i)+(ae+ei+ai-cg-bd-fh) &\neq 0. \\
\lambda &= 2: \hspace{20px} 3(a+e+i)+2(ae+ei+ai-cg-bd-fh) &\neq 0. \\
\lambda &= 3: \hspace{20px} -2(a+e+i)+3(ae+ei+ai-cg-bd-fh)-2 &\neq 0. \\
\lambda &= 4: \hspace{20px} -2(a+e+i)-3(ae+ei+ai-cg-bd-fh) &\neq 0. \\
\lambda &= 5: \hspace{20px} 3(a+e+i)-2(ae+ei+ai-cg-bd-fh)-2 &\neq 0. \\
\lambda &= 6: \hspace{20px} -(a+e+i)-(ae+ei+ai-cg-bd-fh)-2 &\neq 0. \\
\end{align*}

We consider all the possible values of $\tr(M)$ and the coefficient of $\lambda$ to get an exhaustive list of characteristic polynomials: of which the trace $0,\pm1$ cases are as follows. \\ 

First, we let $\mbf{[i,-]}$ denote the set of trace $i$ matrices without eigenvectors in $\SL$. Also, we denote $I$ as the $3 \times 3$ identity matrix in $\SL$.

\subsection{Trace 0 matrices with no eigenvectors in $\SL$}
The characteristic polynomials with no roots in $\F_7$ are
$$\lambda^3+\lambda-1.$$
$$\lambda^3+2\lambda-1.$$
$$\lambda^3+4\lambda-1.$$
Using our notation, for ${k=1, 2, 4}$, we denote the set of matrices with the characteristic polynomial ${\lambda^3+k \lambda -1}$ by $\mbf{[0,k]}$. \\

Note that if we multiply a matrix in $\mbf{[0, -]}$ by $2I$, its trace remains zero and the coefficient of $\lambda$ in the characteristic polynomial will be scaled by a factor of $2^2=4$. This is equivalent to $\lambda$ being replaced by $4 \lambda = \frac{\lambda}{2}$ in the characteristic polynomial (using arithmetic in $\F_7$). Hence, if we select a representative element $M$ of $\mbf{[0, 4]}$, then we will have 

\begin{align*}
    2M &\in \mbf{[0,2]}.\\
    4M = 2 \cdot 2M &\in \mbf{[0,1]}.
\end{align*}

We choose a representative matrix in $\mbf{[0, 4]}$ (for lack of a better criteria) with the maximum number of zeroes possible, 
$$
M = \MatrThree{0}{1}{3}{0}{0}{1}{1}{0}{0}.
$$

This matrix can be checked to be in $\mbf{[0, 4]}$. Note that there are no such representative matrices with exactly $6$ zero entries as these would have exactly one non-zero entry in each column and hence would have an eigenvector. Also, there are no representative matrices with $7$ or more zero entries as these would have a row of all zero entries so would have determinant zero. $M$ has order 57 and $M^{19}=4$. This was confirmed by computer and by hand. This also implies that all matrices in $\mbf{[0, -]}$ have an order of either $19$ or $57=3\times 19$, the 3 resulting from the order of the scaling matrix $2I$. In the next section, one such calculation by hand is given explicitly. \\

The matrix $M$ along with $2M$ and $4M$ can be used to represent 3 conjugacy classes of traceless matrices without any eigenvectors. We will endeavour to prove that these are the only three conjugacy classes.  \\ 

The characteristic equations prove that these are not conjugate. Another way of seeing it is to compare eigenvalues $M$ and $2M$ considered as matrices over a larger field of size $7^3$.




\subsection{Trace 1 matrices without eigenvectors}
One characteristic polynomial with no roots in $\F_7$ will be
$$\lambda^3-\lambda^2-1$$
Using similar notation as earlier, denote the set of matrices associated with this characteristic polynomial as $$\mbf{[1, 0]}$$
(trace is $1$, coefficient of $\lambda$ is $0$). A representative of this set is
$$
M = \MatrThree{0}{1}{0}{0}{1}{-1}{-1}{0}{0}.
$$

This has order $57$, with $M^{19}=4I$. \\

We also have the characteristic polynomial 

$$\lambda^3-\lambda^2+3\lambda-1$$

with no roots in $\SL$. The corresponding set of matrices will be denoted as 
$$\mbf{[1, 3]}.$$
A representative of this set is
$$
M = \MatrThree{0}{1}{4}{0}{1}{5}{1}{0}{0}.
$$

This has order $19$, which was confirmed by hand. \\

Finally, we have the characteristic polynomial 

$$\lambda^3-\lambda^2+5\lambda-1$$

with no roots in $\SL$. The corresponding set of matrices will be denoted as 
$$\mbf{[1, 5]}.$$
A representative of this set is 
$$
M = \MatrThree{0}{3}{2}{0}{1}{1}{1}{0}{0}.
$$
This has order $57$, with $M^{19}=2I$.


\subsection{Trace -1 matrices without eigenvectors}

One characteristic polynomial with no roots in $\SL$ will be

$$\lambda^3+\lambda^2+4\lambda-1$$

The corresponding set of matrices will be 
$$\mbf{[-1, 4]}.$$
A representative of this set is

$$
M = \MatrThree{0}{2}{3}{0}{1}{-1}{1}{0}{0}.
$$

This has order $57$. \\

The second characteristic polynomial with no roots in $\SL$ will be 

$$\lambda^3+\lambda^2+2\lambda-1$$

The corresponding set of matrices will be $\mbf{[-1, 2]}$. A representative of this set is

$$
M = \MatrThree{0}{3}{2}{0}{-1}{-1}{-1}{0}{0}.
$$
The order of $M$ is, which was confirmed by hand and by computer (GC). Also, $M^{19}=2I$, which means that the order is $19$ in $\PSL$.

\subsection{Bijections with other traces}

We get bijections

$$\mbf{[i,-]}\to \mbf{[2i,-]}$$
$$M\mapsto 2M$$

Since multiplying a matrix $M \in \mbf{[i, j]}$ by $2$ causes the the coefficient of $\lambda$ in the corresponding characteristic polynomial to multiply by $2^2 = 4 = \frac{1}{2}$ (in $\F_7$), then the above bijection can be restricted to  

$$\mbf{[i,j]}\to\mbf{[2i,4j]}.$$
This corresponds to 
$$\lambda\mapsto \lambda/2$$
for the characteristic polynomials in $\lambda$. \\

Hence, the sets of matrices $\mbf{[1, -]}, \mbf{[2, -]}, \mbf{[4, -]}$ are very closely related and we have a similar result for the sets of matrices $\mbf{[-1, -]}, \mbf{[-2, -]}, \mbf{[-4, -]}$. \\

Since the matrices in $\mbf{[0, -]}$ can be partitioned into $\mbf{[0,1]}, \mbf{[0, 2]}$ and $\mbf{[0, 4]}$, each non-empty and with distinct corresponding characteristic polynomials, then the three representatives for each of the three subsets belongs to a distinct conjugacy class so there are at least three conjugacy classes that partition $\mbf{[0, -]}$.

This gives at least 18 conjugacy classes.

We will prove in the next section that these 18 conjugacy classes are all of them.

There are some other interesting bijections, which we explore more completely later:
\begin{eg}
$$\mbf{[i, j]}\mapsto \mbf{[j, i]}$$
$$\lambda\mapsto \frac{1}{\lambda}$$
\end{eg}






\section{Orders of matrices in $\mbf{[i, j]}$}

\begin{prop}All matrices in $\SL$ without eigenvectors have order 19 or 57
\end{prop}
\begin{proof}
We summarise the results given in later sections that prove this claim. Let $M\in \mbf{[i, j]}$ have characteristic polynomial $P\in \F_7[x]$. It can be shown that $P$ has three distinct roots in the extended field $\F_{7^3}$ (none of which are in $\F_7$) of the form $\alpha$, $\alpha^7$ and $\alpha^{49}$ so $M$ is diagonalisable over $\F_{7^3}$ with the three roots being the diagonal entries of the corresponding diagonal matrix (the eigenvalues). The determinant of this diagonal matrix must be $\alpha^{57}=1$. Hence, the order of $\alpha$ (and also $M$) is a factor of $57$. Since $\alpha \not\in \F_7$, then the order cannot be $1$ or $3$, so must be $19$ or $57$. 
\end{proof}

Knowing all have order 19 or 57 is sufficient to establish the conjugacy classes in the next section but we would not know which had order 19 and which had order 57 (though we would know that for every one of order 19 there are 2 of order 57).

There are many ways that the orders can be obtained.

To find out the order for each set, it suffices to know it for one matrix in it. Here is one such matrix:

Table of first 19 powers of a matrix with characteristic polynomial $\lambda^{3}+2\lambda -1$.
\ \\
\begin{tabular}{c|c|c|c}\label{powertable}
power&matrix&trace&class\\\hline
1&$
\MatrThree{0}{2}{-1}{0}{0}{2}{2}{0}{0}
$&0&$\mbf{[0,2]}$\\\hline
2&$
\MatrThree{-2}{0}{-3}{-3}{0}{0}{0}{-3}{-2}
$&3&$\mbf{[3,4]}$\\\hline
3&$
\MatrThree{1}{3}{2}{0}{1}{3}{3}{0}{1}
$&3&$\mbf{[3,4]}$\\\hline
4&$
\MatrThree{-3}{2}{-2}{-1}{0}{2}{2}{-1}{-3}
$&1&$\mbf{[1,3]}$\\\hline
5&$
\MatrThree{3}{1}{0}{-3}{-2}{1}{1}{-3}{3}
$&4&$\mbf{[4,3]}$\\\hline
6&$
\MatrThree{0}{-1}{-1}{2}{1}{-1}{-1}{2}{0}
$&1&$\mbf{[1,3]}$\\\hline
7&$
\MatrThree{-2}{0}{-2}{-2}{-3}{0}{0}{-2}{-2}
$&0&$\mbf{[0,2]}$\\\hline
8&$
\MatrThree{3}{3}{2}{0}{3}{3}{3}{0}{3}
$&2&$\mbf{[2,0]}$\\\hline
9&$
\MatrThree{-3}{-1}{3}{-1}{0}{-1}{-1}{-1}{-3}
$&1&$\mbf{[1,3]}$\\\hline
10&$
\MatrThree{-1}{1}{1}{-2}{-2}{1}{1}{-2}{1}
$&3&$\mbf{[3,1]}$\\\hline\end{tabular}
\quad
\begin{tabular}{c|c|c|c}
power&matrix&trace&class\\\hline
11&$
\MatrThree{2}{-2}{3}{2}{3}{-2}{-2}{2}{2}
$&0&$\mbf{[0,2]}$\\\hline
12&$
\MatrThree{-1}{-3}{1}{3}{-3}{-3}{-3}{3}{-1}
$&2&$\mbf{[2,0]}$\\\hline
13&$
\MatrThree{2}{-2}{2}{1}{-1}{-2}{-2}{1}{2}
$&3&$\mbf{[3,1]}$\\\hline
14&$
\MatrThree{-3}{-3}{1}{3}{2}{-3}{-3}{3}{-3}
$&3&$\mbf{[3,4]}$\\\hline
15&$
\MatrThree{2}{1}{-3}{1}{-1}{1}{1}{1}{2}
$&3&$\mbf{[3,1]}$\\\hline
16&$
\MatrThree{1}{-3}{0}{2}{2}{-3}{-3}{2}{1}
$&4&$\mbf{[4,3]}$\\\hline
17&$
\MatrThree{0}{2}{0}{1}{-3}{2}{2}{1}{0}
$&4&$\mbf{[4,3]}$\\\hline
18&$
\MatrThree{0}{0}{-3}{-3}{2}{0}{0}{-3}{0}
$&2&$\mbf{[2,0]}$\\\hline
19&$
\MatrThree{1}{0}{0}{0}{1}{0}{0}{0}{1}
$&3&$\mbf{[3,3]}$\\\hline
20&$
\MatrThree{0}{2}{-1}{0}{0}{2}{2}{0}{0}
$&0&$\mbf{[0,2]}$\\\hline
\end{tabular}

The table of powers of a matrix can also be used to get the orders of all the other sets because a representative from each set is present.

Another way is to show any $\lambda$ satisfying the characteristic equation will also satisfy $\lambda^{19}=1$ as follows:

\begin{prop}\label{02nt}
All matrices in $\mbf{[0, 2]}$ have order 19. 
\end{prop}
\begin{proof}  
We apply some results given later in later sections. Since all matrices $M$ in $\mbf{[0,2]}$ have no eigenvalues in $\F_7$, then by Lemma \ref{EigenVectorWhenExtended}, $M$ has three distinct eigenvalues in $\F_{7^3}$ (none of which are in $\F_7$) and each eigenvalue is a root of the equation $\lambda^3 + 2\lambda - 1 = 0$. \\

By Lemma $\ref{OrderNIffExtSame}$ it suffices to prove that any root $\lambda$ of the above equation has order $19$. We will prove that $\lambda^{21} = \lambda^2$, which is equivalent to $\lambda^{19} = 1$. Below we will use the property that $(a+b)^7 = a^7+b^7$ for all $a, \ b \in \F_{7^3}$. \\

First, we note from the cubic equation that $\lambda^3 = 1 - 2\lambda$. Now we compute $\lambda^{21}$ below.

. We show $\lambda^{21}=\lambda^2$:
$$\begin{aligned}\lambda^{21}&=(1-2\lambda)^7)\\
                            &= 1-2^7\lambda^7\\
                            &= 1 - 2^7 \lambda (\lambda^3)^2\\
                            &=1-2\lambda(1-2\lambda)^2\\
                        &=1-2\lambda(1-4\lambda+4\lambda^2)\\
                        &=1-2\lambda+\lambda^2-\lambda^3\\
                        &=1-2\lambda+\lambda^2-(1-2\lambda)\\
                        &=\lambda^2.
                            \end{aligned}$$
                            
Hence, the order of all three eigenvalues is a factor of $19$. Since $19$ is prime and none of the eigenvalues are equal to $1$ (since the eigenvalues are not in $\F_7$), then the eigenvalues have order $19$.
\end{proof}

This above calculation with the $\lambda$ relation can account for all the matrices in the set at once.\footnote{Other such calculations are included in the Explanatory notes Section.}

We can also use a combination of these strategies, as is demonstrated in the final proof.

\begin{thm}\label{order} All matrices in $\SL$ without eigenvectors have order 19 or 57. Moreover,    
\begin{itemize}
   \item All the elements in each of the sets $\mbf{[0, 2]}$, $\mbf{[1, 3]}$, $\mbf{[2, 0]}$, $\mbf{[3, 1]}$, $\mbf{[3, 4]}$ and $\mbf{[4, 3]}$ have order $19$.
    \item All the elements in each of the sets $\mbf{[0, 1]}$, $\mbf{[0, 4]}$, $\mbf{[1, 0]}$, $\mbf{[1, 5]}$, $\mbf{[2, 5]}$, $\mbf{[2, 6]}$, $\mbf{[4, 0]}$, $\mbf{[4, 6]}$, $\mbf{[5, 1]}$, $\mbf{[5, 2]}$, $\mbf{[6, 2]}$ and $\mbf{[6, 4]}$ have order $57$.
\end{itemize}
\end{thm}
\begin{proof}
We summarise the proof given in Lemma \ref{CharacteriseAllOrder19And57Matrices}. \\

Consider the matrix 

$$
M = \MatrThree{0}{2}{-1}{0}{0}{2}{2}{0}{0}.
$$

We can use Table \ref{powertable} to select powers of $M$ that are in each of the sets $\mbf{[0, 2]}$, $\mbf{[1, 3]}$, $\mbf{[2, 0]}$, $\mbf{[3, 1]}$, $\mbf{[3, 4]}$ and $\mbf{[4, 3]}$, and these matrices will also have order $19$. \\

Scaling these matrices by $2$ or $4$ and using Lemma \ref{OrderOfDoubleAndQuadrupleIs57} gives us matrices in each of the other $12$ sets having order $57$. The theorem then follows by Lemma \ref{OneOrderNThenAll}.

%
\end{proof}

\section{Conjugacy classes}

\begin{center}
    \begin{tabular}{|c|c|c|}
    \hline
        Trace & matrices in $\SL$ without eigenvectors& classes\\
    \hline
    \hline
         0 & $296352 = 3(7^3 \times 3^2 \times 2^5)$&$\mbf{[0, 1]},\mbf{[0, 2]},\mbf{[0, 4]}$\\
    \hline
        1 & $296352 = 3(7^3 \times 3^2 \times 2^5)$&$\mbf{[1, 0]},\mbf{[1, 3]}, \mbf{[1, 5]}$\\
    \hline
        2 & $296352 = 3(7^3 \times 3^2 \times 2^5)$&$\mbf{[2, 0]},\mbf{[2,5]},\mbf{[2,6]}$\\
    \hline
        3 & $197568 = 2(7^3 \times 3^2 \times 2^5)$&$\mbf{[3,4]},\mbf{[3,1]}$\\
    \hline
        4 & $296352 = 3(7^3 \times 3^2 \times 2^5)$&$\mbf{[4,0]},\mbf{[4,6]},\mbf{[4,3]}$\\
    \hline
        5 & $197568 = 2(7^3 \times 3^2 \times 2^5)$&$\mbf{[5,1]},\mbf{[5,2]}$\\
    \hline 
        6 & $197568 = 2(7^3 \times 3^2 \times 2^5)$&$\mbf{[6,2]},\mbf{[6,4]}$\\
    \hline
    \end{tabular}
\end{center}
This was first done by computer and similar results appear online. The text files generated by computer show all the matrices which satisfy the conditions. 
\newline
We proceed to explain and prove these results independent of computational power.

\begin{lem} Let
$$
M = \MatrThree{0}{1}{3}{0}{0}{1}{1}{0}{0}.
$$

The only matrices which commute with $M$ are the $57$ powers of $M$, consequently the size of its conjugacy class is 
$$2^5\times 3^2\times 7^3.$$
\end{lem}
\begin{proof}
The stabiliser of the matrix $M$ can be used to calculate the size of the conjugacy class. Obviously, the stabilisers of $2M$ and $4M$ are identical. For $SM=MS$ we get

$$S=\begin{pmatrix} 
a&b&3b+d\\
d&a-3d&b\\
b&d&a\\
\end{pmatrix}
$$
The determinant 1 condition is
$$a^3+b^3+d^3-3a^2d-3ab^2+2db^2-d^2b-3bad=1$$

One way to see this equation has $57$ solutions is by case bashing (we can reduce the number of cases from $7^3$ to $a=0,1,-1$ since we can scale by $2$ a similar calculation is summarised in the Explanatory notes Section). These must then be the $57$ powers of $M$.

A more indirect approach is as follows. Since these form a subgroup, the number must be divisible by $57$ and also a divisor of $2^5\times3^3\times7^3\times19$. It cannot have a factor of $7$ as $7\times57>7^3$. For $3$, any order 3 matrix not $2I,4I$ must have a one dimensional eigenspace preserved by $M$, a contradiction, while if a matrix has order $9$, its cube has order 3 but cannot equal $I$ since it was order 9, it cannot be $2I$ or $4I$ since $2$ and $4$ are not cubes modulo 7. Of course a matrix with order 27 would cube to a matrix with order 9. Now it suffices to show that no order 2 matrix $S$ commutes with $M$.


Any such matrix must have all eigenvalues $1$ and a two dimensional 1-eigenspace, as otherwise $M$ will have an eigenvector. We show that no such order two matrix exists. Up to conjugacy we can assume
$S$ is of the form
$$
\MatrThree{1}{0}{x}{0}{1}{y}{0}{0}{1}
$$

Consider the last column of $S^2$. The first row entry is $2x$, and the second row entry is $2y$. For $S$ to have order $2$, it follows that $x$ and $y$ are both 0, and so the matrix is the identity matrix, which has order $1$. Therefore there are no matrices that commute with $M$ and have order $2$.

Hence by the Orbit-Stabiliser theorem, the size of the conjugacy classes for $M$, $2M$ and $4M$ are all
$$2^5\times 3^2\times 7^3.$$
\end{proof}

Recall Sylow's theorems from group theory:

\begin{defn}
    \textbf{Sylow $p$-subgroups: } For a prime number $p$, a Sylow $p$-subgroup of a group $G$ is a maximal $p$-subgroup of $G$ (its order is a power of $p$) and is not a proper subgroup of any other $p$-subgroup of $G$.  
\end{defn}
\begin{thm}\label{SylowTheorems}
    \textbf{Sylow's Theorems:} 
    
    \begin{enumerate}
        \item For every prime factor $p$ with multiplicity $n$ of the order of a finite group $G$, there exists a Sylow-$p$ subgroup of $G$ of order $p^n$.
        \item Given a finite group $G$ and a prime number $p$, all Sylow $p$-subgroups of $G$ are conjugate to each other. That is, if $H$ and $K$ are Sylow $p$-subgroups of $G$, then there exists an element $g \in G$ with $g^{-1}Hg=K$.
        \item Let $p$ be a prime factor with multiplicity $n$ of the order of a finite group $G$, so that the order of $G$ can be written as $p^n m$, where $n > 0$ and $p$ does not divide $m$. Let $n_p$ be the number of Sylow $p$-subgroups of $G$. Then the following hold:
        \begin{itemize}
            \item $n_p$ divides $m$, which is the index of the Sylow $p$-subgroup in $G$.
            \item $\ModCmp{n_p}{1}{p}$.
            \item $n_p = |G:N_G(P)|,$ where $P$ is any Sylow $p$-subgroup of $G$ and $N_G$ denotes the normalizer. 
        \end{itemize}
    \end{enumerate}
\end{thm}

\begin{lem}\label{PossibleSubgroupSize}
 If $n_{19}$ is 1 modulo 19 and $n_{19}\mid 2^53^37^3$, then 
 $$n_{19} \in \{1, \  2^5 \cdot 3, \ 7^3, \ 2^4 \cdot 3^2 \cdot 7, \ 2^3 \cdot 3^3 \cdot 7^2, \ 2^5 \cdot 3 \cdot 7^3       \}.$$
\end{lem} 

A further argument will narrow this down to just $2^5 \cdot 3 \cdot 7^3$ in the proof of Theorem \ref{ExactSizeOfConjClass}. This will mean that the normaliser has size $3^2\cdot 19$.


\begin{lem}\label{samestab}
The size of all 18 conjugacy classes above is $2^5 \times3^2\times 7^3$
\end{lem}
\begin{proof}
By the orbit stabiliser theorem it suffices to show the stabilisers of any representative of each conjugacy class have the same size. Take a matrix $M$ of order $19$, the subgroup it generates is a Sylow 19-subgroup and conjugate to all other order 19 subgroups by Sylow's theorem. Therefore for any other conjugacy class, there exists $j$ such that $M^j$ is in that class (or a multiple of something in that class by 2 or 4). Therefore anything that stabilises $M$ will stabilise $M^j$, that is $XMX^{-1}=M$ implies
$XM^jX^{-1}=M^j$. Conversely since 19 is prime, $M^j$ is also a generator and hence $M$ is also a power of $M^j$. Clearly scaling by $2$ does not affect the size of the stabiliser as $2$ commutes with everything.
\end{proof}

\begin{thm}\label{ExactSizeOfConjClass}There are exactly 18 conjugacy classes of matrices without eigenvectors in $\SL$, each have the same size 
$$2^5\times3^2\times7^3$$
which gives a total of
$$2^6\times 3^4\times 7^3.$$
These $18$ conjugacy classes are the sets $\mbf{[0, 1]}$, $\mbf{[0, 2]}$, $\mbf{[0, 4]}$, $\mbf{[1, 0]}$, $\mbf{[1, 3]}$, $\mbf{[1, 5]}$, $\mbf{[2, 0]}$, $\mbf{[2, 5]}$, $\mbf{[2, 6]}$, $\mbf{[3, 1]}$, $\mbf{[3, 4]}$, $\mbf{[4, 0]}$, $\mbf{[4, 3]}$, $\mbf{[4, 6]}$, $\mbf{[5, 1]}$, $\mbf{[5, 2]}$, $\mbf{[6, 2]}$ and $\mbf{[6, 4]}$.
\end{thm}
\begin{proof}
By Lemma \ref{samestab} it suffices to show that we cannot exceed $2^5\times 3^2\times7^3\times18$ matrices without eigenvectors which have orders $19$ or $57$. By Sylow's theorems, there are 5 possible values for the number of order $19$ subgroups (Corollary \ref{PossibleSubgroupSize}), the greatest of which is $2^5\times 3\times 7^3$, these all consist of $18$ distinct matrices without eigenvectors and gives $2^5\times 3\times 7^3\times 18$ such matrices of order 19. Scaling these matrices by 2 and 4 we get the required $2^5\times 3^2\times 7^3\times 18$ matrices, while less would be insufficient.

Now by Proposition \ref{02nt} and similar calculations we see that every matrix without eigenvectors have order 19 or 57.
\end{proof}

\section{Relationships between classes}

Let us analyse this Table \ref{powertable} to see what it tells us.
One thing we notice is that all of the conjugacy classes with order 19 are represented here (this would be an alternate proof that they have order 19 as well). The powers needed to obtain each conjugacy class does not depend on the matrix used only on the conjugacy class. The powers needed for any other class can be deduced from those in one table. The order 19 class is determined by the trace except for when the trace is 3.


\subsection{Bijections between conjugacy classes}

We summarise the bijections resulting from the different powers.

For three generators modulo 19, we get upon repeated application, the bijections cycling between
$$\mbf{[3,4]}\to\mbf{[1,3]}\to\mbf{[2,0]}\to\mbf{[4,3]}\to\mbf{[3,1]}\to\mbf{[0,2]}$$
$$M\mapsto M^2$$
$$M\mapsto M^3$$
$$M\mapsto M^{14}$$
The inverses of these generators, 
$$M\mapsto M^{10}$$
$$M\mapsto M^{13}$$
$$M\mapsto M^{15}$$
traverse these sets in the opposite direction.

Since $7=2^6$, we get for all valid $i,j$, bijections
$$\mbf{[i,j]}\to \mbf{[i,j]}$$
$$M\mapsto M^7$$
and its inverse
$$M\mapsto M^{11}=M^{7^2}.$$

9,4,6 have two three cycles in forward direction, their inverses 17, 5, 16 in the reverse direction
$$\mbf{[3,4]}\to\mbf{[2,0]}\to\mbf{[3,1]}$$
$$\mbf{[1,3]}\to\mbf{[4,3]}\to\mbf{[0,2]}$$

8 have three pairs with inverse 12, as well as 18 which is its own inverse
$$\mbf{[i,j]}\to\mbf{[j,i]}$$

The 18 recovers the bijection in Example \ref{}.

One way in which these bijections are nice is that they are order preserving and commute with the conjugation action (isomorphism of $G$-sets).

\subsection{Simultaneous conjugacy}
Let us make some observations here about the problem of simultaneous conjugacy for commuting matrices.

Is it possible to get a set of commuting matrices where each represents a different conjugacy class?

When are collections of commuting matrices simultaneously conjugate?

If we take an $n$-tuple of commuting elements in $\SL$ $(A_1,\dots , A_n)$ where $A_k\in \mbf{[i_k,j_k]}$ for some $i_k,j_k$ and there is another commuting $n$-tuple $\SL$ $(B_1,\dots , B_n)$ where $B_k\in \mbf{[i_k,j_k]}$, they might not be simultaneously conjugate because the powers might not match. However if the powers match. then they are simultaneously conjugate.

\begin{thm}
Every collection of commuting matrices excluding $2I,4I$ either consists of matrices all with eigenvectors or they consist of matrices all without eigenvectors and all powers of some $M$.
Two such collections are simultaneously conjugate if and only if there exists $M_1$ for the first collection to be powers of $M_1$ and $M_2$ of the same order for the second collection in which all the corresponding powers agree modulo their common order.
\end{thm}
\begin{proof}
If there is a matrix without any eigenvectors in the collection then by Lemma \ref{samestab}, it can only commute with the $57$ powers of some $M$ without eigenvectors. So among the commuting matrices there is a pair $A,B$ for which $B=A^j$ and this relationship is invariant under conjugacy and therefore must hold for any pair of matrices simultaneously conjugate to $A,B$.
\end{proof}


\section{Understanding the subgroup structure}
A maximal subgroup of $\SL$ consists of all block upper triangular matrices of the form
$$
\MatrThree{a}{b}{c}{0}{e}{f}{0}{h}{i}
$$
with determinant 1 and entries in $\F_7$. The number of these is
$$(7^2-1)(7^2-7)7^2=2^5\times 3^2\times 7^3$$

Clearly, for any block upper triangular matrix, the vector 
$$
e_1 := \vecthree{1}{0}{0}
$$
is an eigenvector with eigenvalue $a$.

Any matrix with an eigenvector is conjugate to such a block upper triangular matrix and thus belongs to a subgroup conjugate to the block upper triangular matrices. 

These maximal subgroups intersect  non trivially either in pairs or 3 at a time but not more. These are those with at least 2 or 3 dimensions worth of eigenspaces are are conjugate to 

$$
\MatrThree{a}{0}{c}{0}{e}{f}{0}{0}{i}
$$

or 
$$
\MatrThree{a}{0}{0}{0}{e}{0}{0}{0}{i}
$$

are diagonalisable.

To see what other sorts of subgroups there are, we first search for matrices in $\SL$ that do not have any eigenvectors. We see there are many order 19 and 57 subgroups which contain them and don't intersect any of the other maximal subgroups. So there is another maximal subgroup not of the kind above. However if one extends to $\F_{7^3}$ it is contained in one of the kind above but with entries in a bigger field.

\begin{lem}
$\SL$ does not have any matrices of order $27$.
\end{lem}
\begin{proof}
Let $x\in\SL$ and Consider the extension of $\F$ so that $x$ is diagonalisable. The orders of the eigenvalues must be a factor of $7-1$, $7^2-1$ and $7^3-1$ but these have at most a factor of $3^2$ not $3^3$.
\end{proof}
We definitely know there are elements of order $3$.
The previous lemma also suggests there might be elements of order $3^2$ though this is not the case.

\begin{lem}
There are no elements of order 9 in $\SL$.
\end{lem}
\begin{proof}
One approach is to show that no eigenvalues have order 9 in the extended field of size $7^2$. It is not cubed because order 9 matrices have an eigenvector.
\end{proof}

So far we know from Sylow's theorem that the normaliser $N_{\SL}(P)$ of a subgroup ($P$ of size 19) has size $3^2\times19$, which must contain $P$ together with the center of $\SL$ ($3\times 19$ elements). This means that the size of the normaliser in an subgroup $H$ containing $P$ will have size $19$, $3\times 19$ or $3^2\times19$ depending on whether it is divisible by 3 or 9. 

Questions: The size of normaliser is divisible by $9$, are there any elements of order 9 present?
What is the group of size 9 by 19 up to isomorphism? can we exhibit an isomorphism?
(I know what it should be up to isomorphism)

The other possible values for the number of Sylow $19$-subgroups 
from Lemma \ref{PossibleSubgroupSize} then give us the possible sizes of subgroups that contain $P$, once mutliplied by the size of the normaliser (orbit stabiliser theorem, with $N_{\SL}(P)$ the stabiliser). That is, possible maximal subgroup sizes are 

$2^5\times3\times19$,$2^5\times3^2\times19$,$2^5\times3^3\times19$, $7^3\times 19$, $7^3\times 3\times 19$, $7^3\times 3^2\times 19$ and $2^4 \times3^2 \times7\times 19$,$2^4 \times3^3 \times7\times 19$

\begin{tabular}{c|c|c|c|c|c|c|c}
Lemma\ref{PossibleSubgroupSize}  &$1$&$2^5\times3$&$7^3$& $ 2^4 \cdot 3^2 \cdot 7$&$ 2^3 \cdot 3^3 \cdot 7^2$& $ 2^5 \cdot 3 \cdot 7^3$& \\\hline
$\times 19$&$19$&$2^5\times3\times19$&$7^3\times19$& $ 2^4 \cdot 3^2 \cdot 7\times19$&$ 2^3 \cdot 3^3 \cdot 7^2\times19$& $ 2^5 \cdot 3 \cdot 7^3\times19$&\\
times 3 times 19&$3\times19$&$2^5\times3^2\times19$&$3\times7^3\times19$& $ 2^4 \cdot 3^3 \cdot 7\times19$&$ 2^3 \cdot 3^4 \cdot 7^2\times19$& $ 2^5 \cdot 3^2 \cdot 7^3\times19$&\\
$\times 3^2 \times 19$ &$3^2\times19$&$2^5\times3^3\times19$&$3^2\times7^3\times19$& $ 2^4 \cdot 3^4 \cdot 7\times19$&$ 2^3 \cdot 3^5 \cdot 7^2\times19$& $ 2^5 \cdot 3^3 \cdot 7^3\times19$&\\\hline
\end{tabular}
Removing some obvious inconsistencies gives

\begin{tabular}{c|c|c|c|c|c|c|c}
$\times 19$&$19$&$2^5\times3\times19$&$7^3\times19$& $ 2^4 \cdot 3^2 \cdot 7\times19$&$ 2^3 \cdot 3^3 \cdot 7^2\times19$& $ 2^5 \cdot 3 \cdot 7^3\times19$&\\
$\times 3 \times 19$&$3\times19$&$2^5\times3^2\times19$&$3\times7^3\times19$& $ 2^4 \cdot 3^3 \cdot 7\times19$&$3^4$ too big& $ 2^5 \cdot 3^2 \cdot 7^3\times19$&\\
$\times 3^2 \times 19$ &$3^2\times19$&$2^5\times3^3\times19$&$3^2\times7^3\times19$& $3^4$ too big&$3^5$ too big & whole group&\\
\end{tabular}

{\bf mini questions:} is there a subgroup of index 3? Or index 9?

\begin{lem} If $H$ is a subgroup containing $P$, then the size of $H$ is one of (the entries of the table above, look there first)
$$19,  2^5\cdot 3\cdot 19, 7^3\cdot 19, 2^4\cdot 3^2\cdot 7\cdot 19, 2^3\cdot 3^3\cdot 7^2\cdot 19,2^5\cdot 3\cdot7^3\cdot19,$$
if stabiliser is 19,

$$3\cdot 19,2^5\times3^2\times19,7^3\times 3\times 19, 2^4 \times3^3 \times7\times 19,$$
if stabiliser is $3\times19$

[next one needs to be multiplied by 9]
$$19,  2^5\cdot 3\cdot 19, 7^3\cdot 19, 2^4\cdot 3^2\cdot 7\cdot 19, 2^3\cdot 3^3\cdot 7^2\cdot 19,2^5\cdot 3\cdot7^3\cdot19,$$

\end{lem}
\begin{proof}
The size of $H$ divides the size $\SL$. The number of Sylow-19 subgroups of $H$ is 1 modulo 19 and divides $\SL$, by Lemma \ref{PossibleSubgroupSize} must be one of 
 $$ \{1, \  2^5 \cdot 3, \ 7^3, \ 2^4 \cdot 3^2 \cdot 7, \ 2^3 \cdot 3^3 \cdot 7^2, \ 2^5 \cdot 3 \cdot 7^3       \}.$$
To get the size of $H$ from this we use the orbit stabiliser theorem for the conjugation action on the Sylow-19 subgroups of $H$. The stabiliser is the normaliser of $P$ in $H$, which has size $19$ or $57$. There is a single orbit by Sylow's theorems.
\end{proof}

\section{Consequences for $\PSL$}

The 18 conjugacy classes in $\SL$ collapse to 6 conjugacy classes once the scaling by 2 and 4 are identified.

We denote these classes using square notation on the class with order 19

$$\mbf{[3,4]}, \mbf{[1,3]},\mbf{[2,0]},\mbf{[4,3]},\mbf{[3,1]},\mbf{[0,2]}$$

The remaining 6  classes are still related by the same bijections before (state them again)

\begin{thm}
$$\mbf{[3,4]}\to\mbf{[1,3]}\to\mbf{[2,0]}\to\mbf{[4,3]}\to\mbf{[3,1]}\to\mbf{[0,2]}$$
$$M\mapsto M^2$$
$$M\mapsto M^3$$
$$M\mapsto M^{14}$$
The inverses of these generators, 
$$M\mapsto M^{10}$$
$$M\mapsto M^{13}$$
$$M\mapsto M^{15}$$
traverse these sets in the opposite direction.

Since $7=2^6$, we get for all valid $i,j$, bijections
$$\mbf{[i,j]}\to \mbf{[i,j]}$$
$$M\mapsto M^7$$
and its inverse
$$M\mapsto M^{11}=M^{7^2}.$$

9,4,6 have two three cycles in forward direction, their inverses 17, 5, 16 in the reverse direction
$$\mbf{[3,4]}\to\mbf{[2,0]}\to\mbf{[3,1]}$$
$$\mbf{[1,3]}\to\mbf{[4,3]}\to\mbf{[0,2]}$$

8 have three pairs with inverse 12, as well as 18 which is its own inverse
$$\mbf{[i,j]}\to\mbf{[j,i]}$$
\end{thm}

\newpage
\section{Expository and explanatory notes}
We include in this section omitted details, further explanations and supplementary arguments in parallel with the main article.

\subsection{Group Sizes}
\begin{lem}\label{first_col_SL_anything}
Let $\vv{v} \neq \vv{0}$ be a vector in $\F_7^n$ where $n \geq 2$ is a positive integer. Then there exists a matrix $M \in \GenSL$ with its first column being $\vv{v}$.
\end{lem}
\begin{proof}
We can extend $\vv{v}$ to the basis $\{ \vv{v_1}, \vv{v_2}, \ldots, \vv{v_n} \}$ of $\F_7^n$ where $\vv{v} = \vv{v_1}$. Then the matrix $A$ with columns $\vv{v_1}, \vv{v_2}, \ldots, \vv{v_n}$ going left to right is invertible and suppose that $\det(A) = k \neq 0$. Thus, if we consider the matrix $M$ which is formed by changing the second column from $\vv{v_2}$ to $\frac{\vv{v_2}}{k}$, then 
$$
\det(M) = \frac{1}{k} \cdot \det(A) = \frac{1}{k} \cdot k = 1.
$$
Hence, we have found the required $M$.
\end{proof}

\begin{lem}\label{GL2F7_size}
The size of $\GLTwo$ is $(7^2-1)(7^2-7)$.
\end{lem}
\begin{proof}
To obtain a matrix $M \in \GLTwo$, the first column $\vv{v_1}$ must be a non-zero vector in $\F_7^2$, giving us $7^2-1$ options. The second column $\vv{v_2}$ must be linearly independent to $\vv{v_1}$. There are $7$ vectors linearly dependent to $\vv{v_1}$, (which are of the form $c \vv{v_1}$ for all $0 \leq c \leq 6$). Therefore, there are $7^2-7$ options for $\vv{v_2}$. Multiplying these numbers of options gives us ${(7^2-1)(7^2-7)}$ invertible matrices $M$, and so this is the size of $\GLTwo$.
\end{proof}

\begin{lem}\label{GL3F7_size}
The size of $\GL$ is $(7^3-1)(7^3-7)(7^3 - 7^2)$.
\end{lem}
\begin{proof}
We use a similar approach as in the proof of Lemma \ref{GL2F7_size}. To obtain a matrix $M \in \GL$, the first column $\vv{v_1}$ must be a non-zero vector in $\F_7^3$, giving us $7^3-1$ options. The second column $\vv{v_2}$ must be linearly independent to $\vv{v_1}$. Since there are $7$ vectors linearly dependent to $\vv{v_1}$, (which are of the form $c \vv{v_1}$ for $c=0, 1 \ldots, 6$), then there are $7^3-7$ options for $\vv{v_2}$. \\

The third column $\vv{v_3}$ must be linearly independent to $\vv{v_1}$ and $\vv{v_2}$. Since $\vv{v_1}$ and $\vv{v_2}$ are themselves independent, then there are $7^2$ vectors linearly dependent to $\vv{v_1}$ and $\vv{v_2}$, (which are of the form $c \vv{v_1} + d \vv{v_2}$ for all $0 \leq c, d \leq 6$). Thus, there are $7^3-7^2$ options for $\vv{v_3}$. Multiplying these numbers of options gives us ${(7^3-1)(7^3-7)(7^3-7^2)}$ invertible matrices $M$, and so this is the size of $\GL$.
\end{proof}

\begin{lem}\label{SL_size_from_GL}
For all positive integers $n \geq 2$ we have $|\GenSL| = \frac{1}{6} \cdot |\GenGL|$.
\end{lem}
\begin{proof}
We can partition $\GenGL$ into equivalence classes each of size $6$ where ${M \sim N}$ if and only if the first column $\vv{v}$ of $M$ is of the form $c \vv{w}$ for some $1 \leq c \leq 6$ where $\vv{w}$ is the first column of $N$. \\

Consider an equivalence class $A$. Note that by definition of our equivalence classes, the determinants of the matrices in $A$ will be of the form $k \times \det(M)$ for all $1 \leq k \leq 6$ and where $M$ is one representative of $A$. Since $M \in \GenGL$, then $\det(M) \neq 0$. Hence, exactly one matrix in $A$ will have determinant one (corresponding to $k = \det(M)^{-1}$). \\

If we apply this to all equivalence classes, we obtain that one sixth of the elements of $|\GenGL|$ are in $|\GenSL|$, so the result follows.
\end{proof}

\begin{lem}\label{SL2_size}
The size of $\SLTwo$ is $(7^2-1) \cdot 7$.
\end{lem}
\begin{proof}
We present two proofs. The first proof is an immediate application of Lemmas \ref{GL2F7_size} and \ref{SL_size_from_GL} where we have 
$$
|\SLTwo| = \frac{1}{6} \cdot |\GLTwo| = \frac{1}{6} \cdot (7^2-1)(7^2-7) = (7^2-1) \cdot 7.
$$
For the second proof, we will apply the orbit-stabiliser theorem. Let $A$ be the set of all non-zero vectors $\vv{v} \in \F_7^2$. Then ${|B| = 7^2-1}$. Let 
$$
\vv{u} = \vectwo{1}{0} \in B.
$$
Note that $B$ is a G-set under the group action of left multiplication by any ${M \in \SLTwo}$. For any ${\vv{w} \in B}$, by Lemma \ref{first_col_SL_anything}, we can choose a matrix ${M \in \SLTwo}$ with its first column equal to $\vv{w}$, so we will have ${A \vv{u} = \vv{w}}$. \\

Thus, the size of the orbit of $\vv{u}$ under the group action is ${|B| = 7^2-1}$. Let $G_{\vv{u}}$ be the set of matrices in $\SLTwo$ which stabilise $\vv{u}$, By the orbit stabiliser theorem we must have ${|\SLTwo| = |G_{\vv{u}}||B|}$. \\

Since for all $A \in \SLTwo$ we have $A \vv{u} = \vv{w}$ where $\vv{w}$ is the first column of $A$, then the matrices in $G_{\vv{u}}$ must be of the form
$$
M = \MatrTwo{1}{j}{0}{k} \in \SLTwo.
$$
Since we have
$$
\det(M) = 1 \cdot k - 0 \cdot j = k,
$$
then $M \in \SLTwo$ if and only if $k = 1$. Then we have $7$ possible choices for $j$, so ${|G_{\vv{u}}| = 7}$. This gives us 
$$
|\SLTwo| = |G_{\vv{u}}||B| = 7 \cdot (7^2 - 1).
$$
\end{proof}

\begin{prop} The size of $\SL$ is
$$\frac{1}{6} \cdot (7^3-1)(7^3-7)(7^3-7^2) = 2^5\times3^3\times 7^3\times 19.$$

\end{prop}
\begin{proof}
Like before, we present two proofs. The first proof is an immediate application of Lemmas \ref{GL3F7_size} and \ref{SL_size_from_GL} where we have
$$
|\SL| = \frac{1}{6} \times |\GL| = \frac{1}{6} \cdot (7^3-1)(7^3-7) (7^3-7^2).
$$
For the second proof, again we will apply the orbit-stabiliser theorem. Let $B$ be the set of all non-zero vectors $\vv{v} \in \F_7^3$. Then ${|B| = 7^3-1}$. Let
$$
\vv{u} = \vecthree{1}{0}{0} \in B.
$$
Note that $B$ is a G-set under the group action of left multiplication by any ${M \in \SL}$. For any ${\vv{w} \in B}$, by Lemma \ref{first_col_SL_anything}, we can choose a matrix ${M \in \SL}$ with its first column equal to $\vv{w}$, so we will have ${A \vv{u} = \vv{w}}$. \\

Thus, the size of the orbit of $\vv{u}$ under the group action is ${|B| = 7^3-1}$. Let $G_{\vv{u}}$ be the set of matrices in $\SL$ which stabilise $\vv{u}$, By the orbit stabiliser theorem we must have ${|\SL| = |G_{\vv{u}}||B|}$. \\

Since for all $A \in \SL$ we have $A \vv{u} = \vv{w}$ where $\vv{w}$ is the first column of $A$, then the matrices in $G_{\vv{u}}$ must be of the form
$$
M = \MatrThree{1}{b}{c}{0}{e}{f}{0}{h}{i}.
$$
Since we have
$$
\det(M) = ei-fh,
$$
then $M \in \SL$ if and only if $ei-fh = 1$, which is equivalent to the sub-matrix
$$
N = \MatrTwo{e}{f}{h}{i}
$$
being in $\SLTwo$. By Lemma \ref{SL2_size}, we have ${|\SLTwo| = (7^2-1) \cdot 7}$ possible choices for ${e, f, h, i}$ and a further $7^2$ choices for $b, c$. Therefore, we have 
$$
|G_{\vv{u}}| = 7^3 \cdot (7^2 - 1).
$$
This finally gives us
$$
|\SL| = |G_{\vv{u}}||B| = 7^3 \cdot (7^2-1)(7^3-1) = 7^2 \cdot (7^3-7)(7^3-1) = \frac{1}{6} \cdot (7^3-1)(7^3-7) (7^3-7^2).
$$
\end{proof}

\subsection{Irreducible characteristic polynomials}

Consider the matrix 
$$
M = \MatrThree{a}{b}{c}{d}{e}{f}{g}{h}{i} \in \SL.
$$
To find the eigenvalues of $M$, we solve the equation $\mathrm{det}(M - \lambda I) = 0$. Since
$$
M = \MatrThree{a-\lambda}{b}{c}{d}{e - \lambda}{f}{g}{h}{i-\lambda}
$$
then the determinant equation becomes 
$$(a-\lambda)(e-\lambda)(i-\lambda) + bfg + cdh - c(e-\lambda)g - fh(a-\lambda) - bd(i-\lambda) = 0.$$ 
This expands as 
$$ aei+bfg + cdh-ceg-afh-bdi + (-1)\lambda^3 + (a+e+i)\lambda^2 + (-ae-ei-ai+cg+fh+bd)\lambda = 0.$$
Also, since
$$
aei + bfg + cdh - ceg - afh - bdi = \det(M) = 1
$$
then the equation rearranges to the characteristic polynomial equation for $M$ of
$$
\lambda^3 - (a+e+i)\lambda^2 + (ae+ei+ai-cg-bd-fh)\lambda -1 = 0.
$$
We can alternatively use the fact ${\tr(M) = a+e+i}$ to write the equation as
$$
\lambda^3 - \tr(M) \lambda^2 + (ae+ei+ai-cg-bd-fh)\lambda -1 = 0.
$$

In order for $M$ to have no eigenvalues, this equation cannot hold for any integer $1 \leq \lambda \leq 6$ (we already know that $\lambda$ cannot be zero because $M$ is invertible). By substituting each value of $\lambda$ into the equation, this gives us the six equations in the introduction that must not hold in order for $M$ to have no eigenvalues.

\subsection{Simple Properties of the Group}
\begin{lem}\label{DoublingAndQuadruplingStillInSL}
Suppose that $M \in \SL$. Then $2M$ and $4M$ are also in $\SL$.
\end{lem}
\begin{proof}
Since $M \in \SL$, then $\det(M) = 1$. We then have
$$
\det(2M) = 2^3 \cdot \det(M) = 8 \cdot 1 = 1
$$
using $\F_7$ arithmetic. Similarly, we have
$$
\det(4M) = 4^3 \cdot \det(M) = 64 \cdot 1 = 1.
$$
Hence, $2M$ and $4M$ are in $\SL$.
\end{proof}

\begin{lem}\label{OrderOfDoubleAndQuadrupleIs57}
Suppose that $M \in \SL$ and $M$ has order $19$. Then $2M$ and $4M$ have order $57$.
\end{lem}
\begin{proof}
We have $M^{19} = I$. Since $M$ has order $19$, then using $\F_7$ arithmetic we have
$$(2M)^{57} = 2^{57} \cdot (M^{19})^3 = (2^3)^{19} \cdot I = I.$$ 
Thus, the order of $2M$ must be a factor of $57$. Note that 
$$(2M)^{19} = 2^{19} \cdot M^{19} = (2^3)^6 \cdot 2 \cdot I = 2I \neq I.$$
Also, 
$$(2M)^3 = 2^3 \cdot M^3 = M^3 \neq I.$$ 
because $M$ has order $19$. Hence, $2M$ has order $57$. We perform similar computations for the matrix $4M$. We have
$$(4M)^{57} = 2^{2 \cdot 57} \cdot (M^{19})^3 = (2^3)^{38} \cdot I = I.$$ 
This means that the order of $4M$ must be a factor of $57$. Note that 
$$(4M)^{19} = 2^{38} \cdot M^{19} = (2^3)^{12} \cdot 4 \cdot I = 4I \neq I.$$ 
Also, 
$$(4M)^3 = (2^3)^2 \cdot M^3 = M^3 \neq I$$ 
because $M$ has order $19$. Hence, $4M$ has order $57$. 
\end{proof}

\begin{lem}\label{2MAnd4MNoEigenvaluesIfMDoesNotHave}
Suppose that $M \in \SL$ and $M$ has no eigenvalues in $\F_7$. Then $2M$ and $4M$ have no eigenvalues in $\F_7$.
\end{lem}
\begin{proof}
Let $t = 2$ or $4$. Suppose for the sake of contradiction that $tM \vv{v} = \lambda \vv{v}$ for some eigenvector $\vv{v} \neq \vv{0}$ in $\F_7^3$ and some eigenvalue $\lambda \in \F_7$. Then $M \vv{v} = \frac{\lambda}{t} \vv{v}$ so $M$ has eigenvalue $\frac{\lambda}{t} \in \F_7$, which is a contradiction.
\end{proof}

\subsection{Preliminary Results on sets of the form [i, j]}
\begin{lem}\label{M2MBijection}
Suppose that $A = \mbf{[i, j]} \subseteq \SL$ and no matrices in $A$ have any eigenvalues in $\F_7$. Then the set $B = \mbf{[2i, 4j]} \subseteq \SL$ has the same size as $A$ and no matrices in $B$ have any eigenvalues in $\F_7$. Also, there is a bijection from $A$ to $B$ defined by $M \mapsto 2M$.
\end{lem}
\begin{proof}
Let 
$$C = \{ 2M \mid M \in A \} \subseteq \SL. $$
Consider any $M \in A$ where
$$
M = \MatrThree{a}{b}{c}{d}{e}{f}{g}{h}{i}.
$$ 
Let $p(\lambda)$ and $q(\lambda)$ be the characteristic polynomials of $M$ and $2M$ respectively. Since the coefficient of $\lambda^2$ in $p$ is $-\tr(M)$, then similarly the coefficient of $\lambda^2$ in $q$ is $-\tr(2M) = -2 \cdot \tr(M)$. The coefficient of $\lambda$ in $p$ is $ae+ie+ai-cg-bd-fh$. Note that in $2M$, each variable in the previous expression doubles so the coefficient of $\lambda$ in $q$ is four times that of $p$. Thus, $2M \in B$, so $C \subseteq B$. We also get a bijective mapping from $A$ to $C$ defined by $M \mapsto 2M$. \\

Let 
$$
D = \{ 4M \mid M \in B \} \subseteq \SL.
$$
Consider any $M \in B$. Define $p(\lambda)$ as earlier and let $r(\lambda)$ be the characteristic polynomial of $4M$. Using similar reasoning as before, the coefficient of $\lambda^2$ in $r$ is four times (or half times in $\F_7$) that of $p$ and the coefficient of $\lambda$ in $r$ is $4^2 = 2 = \frac{1}{4}$ times that of $p$. Hence, $4M \in A$ and so $D \subseteq A$. We also get a bijective mapping from from $B$ to $D$ defined by $M \mapsto 4M$. \\

We have $|A| = |C|$, $|B| = |D|$, $C \subseteq B$ and $D \subseteq A$. Thus, $|C| \leq |B|$ and $|D| \leq |A|$, which forces $|A|=|B|=|C|=|D|$ and so there is a bijection from $A$ to $B$ defined by $M \mapsto 2M$. Using this bijection and  Lemma \ref{2MAnd4MNoEigenvaluesIfMDoesNotHave}, we find that none of the matrices in $B$ have any eigenvalues in $\F_7$.

\end{proof}

\begin{lem}\label{18PairsRepresentAllNonEigenMatrices}
The following $18$ sets represent all the matrices in $\SL$ with no eigenvalues in $\F_7$: $\mbf{[0, 1]}$, $\mbf{[0, 2]}$, $\mbf{[0, 4]}$, $\mbf{[1, 0]}$, $\mbf{[1, 3]}$, $\mbf{[1, 5]}$, $\mbf{[2, 0]}$, $\mbf{[2, 5]}$, $\mbf{[2, 6]}$, $\mbf{[3, 1]}$, $\mbf{[3, 4]}$, $\mbf{[4, 0]}$, $\mbf{[4, 3]}$, $\mbf{[4, 6]}$, $\mbf{[5, 1]}$, $\mbf{[5, 2]}$, $\mbf{[6, 2]}$ and $\mbf{[6, 4]}$.
\end{lem}
\begin{proof}
Recall that $\mbf{[i, j]}$ is defined as the set of all matrices in $\SL$ with the characteristic polynomial $\lambda^3-i\lambda^2+j\lambda-1$ where $i, j \in \F_7$. We found earlier that all the trace $0$ matrices in $\SL$ with no eigenvalues in $\F_7$ (where $i = 0$) are the sets $\mbf{[0,1]}$, $\mbf{[0, 2]}$ and $\mbf{[0, 4]}$. Also, we showed earlier that all such trace $1$ matrices (where $i = 1$) are the sets $\mbf{[1,0]}$, $\mbf{[1,3]}$ and $\mbf{[1, 5]}$. Furthermore, all such trace $-1$ matrices (where $i = -1$) were found to be the sets $\mbf{[-1, 2]}$ and $\mbf{[-1, 4]}$. \\

Note that $\lambda$ is a root of the polynomial $\lambda^3 - i\lambda^2 + j\lambda -1$ over $\F_7$ if and only if $2\lambda$ is a root of the polynomial $\mu^3-2i\mu^2 + 4j\mu -1$ because 
$$(2\lambda)^3-(2i)(2\lambda)^2+(4j)(2\lambda)-1=\lambda^3-i\lambda^2+j\lambda-1$$ 
using arithmetic in $\F_7$. Therefore, the set $\mbf{[i, j]}$ consists only of matrices with no eigenvalues in $\F_7$ if and only if the same is true for $\mbf{[2i, 4j]}$ and for $\mbf{[\frac{i}{2}, \frac{j}{4}]}$ = $\mbf{[4i, 2j]}$. \\

Hence, a set of the form $\mbf{[2, j]}$ having matrices all without eigenvalues in $\F_7$ is equivalent to the same being true for $\mbf{[1, 2j]}$. This forces $j = 0, \ 5$ or $6$, giving us the sets $\mbf{[2,0]}$, $\mbf{[2,5]}$ and $\mbf{[2,6]}$. For sets of the form $\mbf{[3,j]}$ we can equivalently look at $\mbf{[6,4j]}$, which forces $j=1$ or $4$, giving us the sets $\mbf{[3,1]}$ and $\mbf{[3,4]}$. For sets of the form $\mbf{[4,j]}$ we can equivalently look at $\mbf{[1,4j]}$, which forces $j=0, 3$ or $6$, giving us the sets $\mbf{[4,0]}$, $\mbf{[4,3]}$ and $\mbf{[4,6]}$. For sets of the form $\mbf{[5,j]}$ we can equivalently look at $\mbf{[-1, 2j]}$, which forces $j=1$ or $2$, giving us the sets $\mbf{[5,2]}$ and $\mbf{[5,4]}$. \\

Thus, these $18$ sets consist of all the matrices in $\SL$ without eigenvalues in $\F_7$.
\end{proof}



\subsection{Eigenvalues}

\begin{lem}\label{noRepeatRoot}
Let $p$ be a monic cubic polynomial with coefficients in $\mathbb{F}_7$ which is irreducible over $\mathbb{F}_7$. Then $p$ has no repeated roots in $\mathbb{F}_{7^3}$.
\end{lem}
\begin{proof}
Suppose that there exists such a $p$ which has a repeated root $\alpha \in \mathbb{F}_{7^3}$. Since $p$ is irreducible, then it must be the minimal polynomial of $\alpha$. The derivative $p'$ will also have the root $\alpha$. But $\deg(p') < \deg(p)$, which contradicts the minimality of $\deg(p)$.
\end{proof}

\begin{lem}\label{EigenVectorWhenExtended}
Let $M \in \SL$ such that it has no eigenvalues in $\mathbb{F}_7$. Then $M$ must have three distinct eigenvalues in $\mathbb{F}_{7^3}$ which are of the form $\lambda$, \ $\lambda^7$ and $\lambda^{49}$ for some $\lambda \in \F_{7^3}$. 
\end{lem}
\begin{proof}
The characteristic polynomial $p(\lambda)$ of $M$ is monic and cubic in $\lambda$ with no roots in $\F_7$. Then $p$ must have no linear factors with coefficients in $\F_7$ so $p$ is irreducible over $\F_7$. By Lemma \ref{noRepeatRoot}, $p$ has no repeated roots in $\F_{7^3}$. Also, $p$ must split into linear factors over $\F_{7^3}$, so $p$ must have three distinct roots (which are eigenvalues) in $\F_{7^3}$. Since $p$ is monic and irreducible over $\F_7$, its three distinct roots are of the form $\lambda$, $\lambda^7$ and $\lambda^{7^2}=\lambda^{49}$ for some $\lambda \in \F_{7^3}$.
\end{proof}

\begin{lem}\label{ElementLargerFieldOrderFixedCoprime}
Consider a $\lambda \in \F_{7^3}$ with order $n \geq 1$. If $k$ is a positive integer such that $\gcd(k, n)=1$, then $\lambda^k$ also has order $n$.
\end{lem}
\begin{proof}
Let $t$ be the order of $\lambda^k$. Since 
$$(\lambda^k)^n=\lambda^{kn}=(\lambda^n)^k=1^k=1,$$ 
then $t$ divides $n$. Also, since
$$1=(\lambda^k)^t=\lambda^{kt}$$ 
then $kt$ is a multiple of $n$. Since $\gcd(k, n)=1$, then $t$ is a multiple of $n$, so we must have $t=n$.
\end{proof}

\begin{lem}\label{CommutingMatricesEigenspacesPreserved}
Suppose that $A$ and $B$ are commuting matrices in $\SL$ where $A$ has an eigenvalue $\lambda \in \F_7$ and corresponding eigenspace $E_\lambda$ with dimension one. Then $B$ will have an eigenvalue $\mu \in \F_7$ and corresponding eigenspace $E_\mu$ which satisfies $E_\lambda \subseteq E_\mu$.
\end{lem}
\begin{proof}
Suppose that $\vv{v} \neq \vv{0}$ is in $E_\lambda$. Then $A \vv{v}=\lambda \vv{v}$ and so 
$$A(B \vv{v})=(AB) \vv{v}=(BA) \vv{v}=B(A \vv{v})=B(\lambda \vv{v})=\lambda (B \vv{v}).$$ 
Since $B$ is invertible, then $B \vv{v} \neq \mbf{0}$. From above, we obtain $B \vv{v} \in E_\lambda$, so $B \vv{v} = \mu \vv{v}$ for some $\mu \in \F_7$ (since $E_\lambda$ has dimension one). Hence $\mu$ is an eigenvalue of $B$ and $\vv{v} \in E_\mu$. Since $E_\lambda$ has dimension one, then this gives us that $E_\lambda \subseteq E_\mu$.
\end{proof}

\begin{lem}\label{OrderNIffExtSame}
Consider an $M \in \SL$ with no eigenvalues in $\F_7$. Using Lemma \ref{EigenVectorWhenExtended}, consider the three distinct eigenvalues of $M$ in $\F_{7^3}$. Then $M$ has order $n$ if and only if all three of these eigenvalues have order $n$. Also, the order of $M$ must be a factor of $57$.
\end{lem}
\begin{proof}
Let $\lambda_1, \ \lambda_2$ and $\lambda_3$ be the three distinct eigenvalues of $M$ in $\F_{7^3}$. Then, $M$ can be diagonalised in the form $M = XDX^{-1}$ for some $X \in \GLLarge$ and where 
$$
D = \MatrThree{\lambda_1}{0}{0}{0}{\lambda_2}{0}{0}{0}{\lambda_3} \in \GLLarge.
$$

Note that for all positive integers $k$, $M^k = XD^k X^{-1}$ where we have

$$D^k=\begin{pmatrix} 
{\lambda_1}^k&0&0\\
0&{\lambda_2}^k&0\\
0&0&{\lambda_3}^k\\
\end{pmatrix}
$$

Now, $M^k = I$ is equivalent to $XD^k X^{-1} = I$, or $D^k = X^{-1} I X = I$. Also, $D^k=I$ is equivalent to the orders of all three eigenvalues being a factor of $k$. Let the orders of $\lambda_1$, $\lambda_2$ and $\lambda_3$ be $n_1$, $n_2$ and $n_3$ respectively. Then the smallest value of $k$ such that $D^k=I$ will be $\lcm(n_1, n_2, n_3)$, and this will be the order of $M$.

Next, note that

$$\det(D)=\det(X) \cdot \det(D) \cdot \frac{1}{\det(X)} = \det(XDX^{-1})=\det(M)=1.$$ 

By Lemma \ref{EigenVectorWhenExtended}, we obtain that $(\lambda_1, \lambda_2, \lambda_3)$ is a permutation of $(\lambda, \lambda^7, \lambda^{49})$ for some $\lambda \in \F_{7^3}$. Since $\det(D)=1$, then 
$$\lambda \cdot \lambda^7 \cdot \lambda^{49} = \lambda^{57} = 1.$$
Thus, we also have, $(\lambda^7)^{57} = (\lambda^{49})^{57} = 1$. Hence, $n_1$, $n_2$ and $n_3$ are factors of $57$ so the order of $M$ which equals $\lcm(n_1, n_2, n_3)$ will also be a factor of $57$. This proves the second statement. \\

To prove the first statement, it remains to show that ${n_1=n_2=n_3}$ (since then we would have ${\lcm(n_1, n_2, n_3)=n_1=n_2=n_3}$). Since $\lambda$ has an order $m$ that is a factor of $57$, then we have ${\gcd(m, 7) = \gcd(m, 49) = 1}$. From Lemma \ref{ElementLargerFieldOrderFixedCoprime}, the orders of $\lambda^7$ and $\lambda^{49}$ must also be $m$, so we are done.
\end{proof}

\subsection{Matrix Powers and Orders}
\begin{lem}\label{orderPowerCoprimeFixed}
Consider a $M \in \SL$ with order $n \geq 1$. If $k$ is a positive integer such that $\gcd(k,n)=1$ then $M^k$ also has order $n$.
\end{lem}
\begin{proof}
Let $t$ be the order of $M^k$. Since $M^{kn}=(M^k)^n=I^n=I$, then $t$ divides $n$. Also $I=(M^k)^t=M^{kt}$, so $kt$ is a multiple of $n$. Since $\gcd(k,n)=1$, then $t$ is a multiple of $n$, so we must have $t=n$. 
\end{proof}

\begin{lem}\label{EigenvalueScale}
Consider a $M \in \SL$ with eigenvector $\vv{v} \neq \vv{0}$ in $\F_{7^3}^3$ and corresponding eigenvalue $\lambda \in \F_{7^3}$. Then for any positive integer $k$, the matrix $M^k$ has eigenvector $\vv{v}$ and corresponding eigenvalue $\lambda^k$.
\end{lem}
\begin{proof}
We will use induction.  Suppose that $M^t \vv{v} = \lambda^t \vv{v}$ for some positive integer $t$ (this is true for the base case of $t = 1$). Then 
$$M^{t+1} \vv{v} = M(M^t \vv{v}) = M(\lambda^t \vv{v}) = \lambda^t (M \vv{v}) = \lambda^{t+1} \vv{v}.$$ 
Hence by induction, $M^k \vv{v} = \lambda^k \vv{v}$ for all positive integers $k$.
\end{proof}

\begin{lem}\label{EigenCoprime}
Consider a $M \in \SL$ with order $n \geq 2$ and no eigenvalues in $\F_7$. Then for all positive integers $k$ where $1 \leq k \leq n-1$ and $\gcd(k, n) = 1$, the matrix $M^k$ has no eigenvalues in $\F_7$. 
\end{lem}
\begin{proof}
Suppose that $M^k$ has an eigenvector $\vv{v} \neq \vv{0}$ in $\F_7^3$ and corresponding eigenvalue $\lambda \in \F_7$, Then, $M^k v = \lambda v$. Since $\gcd (k, n) = 1$ then there exists a positive integer $t$ such that $\ModCmp{kt}{1}{n}$. Thus, $kt = sn+1$ for some non-negative integer  $s$. By Lemma \ref{EigenvalueScale}, the matrix 
$$(M^k)^t = M^{sn+1} = (M^n)^s \cdot M = M$$ 
has $\vv{v}$ as an eigenvector, contradiction.
\end{proof}

\begin{lem}\label{PrimeNoEigen}
Consider a $M \in \SL$ with order $p$ where $p$ is prime and where $M$ has no eigenvalues in $\F_7$. Then all the matrices of the form $M^k$ for $1 \leq k \leq p-1$ have no eigenvalues in $\F_7$.

\end{lem}
\begin{proof}
Since $p$ is prime then $\gcd(p, k) = 1$ for $1 \leq k \leq p-1$ so the result follows from Lemma \ref{EigenCoprime}.
\end{proof}

\begin{lem}\label{PowerMatrixSetSameSizeCoprimeOrder}
Let $n \geq 2$ be a fixed positive integer and $A$ be a set of matrices in $\SL$, all of which have order $n$. For a fixed positive integer $k$ where $2 \leq k \leq n$ and $\gcd(k, n) = 1$, let 
$$
B = \{ M^k \mid M \in A \}.
$$
Then $|A| = |B|$.
\end{lem}
\begin{proof}
Since $\gcd(k, n) = 1$, there exists a positive integer $t$ such that ${\ModCmp{kt}{1}{n}}$. Thus, $kt = sn+1$ for some non-negative integer $s$. Since the mapping from $A$ to $B$ defined by $M \mapsto M^k$ is surjective by definition, then $|A| \geq |B|$. Consider an arbitrary $N = M^k \in B$ for some $M \in A$. Note that 
$$N^t = (M^k)^t = M^{sn+1} = (M^n)^s \cdot M = M \in A.$$ 
Thus, the mapping from $B$ to $A$ defined by $M \mapsto M^t$ is surjective, so $|A| \leq |B|$. Hence, $|A| = |B|$. 
\end{proof}

\begin{lem}\label{DifferentSubgroupsDisjoint}
Let $G$ and $H$ be different subgroups of $\SL$ with order $19$. Then, $G \cap H = I$.
\end{lem}
\begin{proof}
Since $G$ and $H$ both contain the identity matrix $I$, then $I \in G \cap H$. Assume for the sake of contradiction that for some non-identity matrix $A \in \SL$, we have $A \in G \cap H$. Since the order $19$ of $G$ and $H$ is a prime number, then both subgroups are cyclic and $A$ is a generator of both $G$ and $H$ (where $A$ has order $19$). Hence, we must have 
$$
G = H = \{ A^k \mid 0 \leq k \leq 18 \}
$$
which contradicts $G$ and $H$ being distinct.
\end{proof}

\begin{lem}\label{AllOrder19SubgroupsAreSylow}
The Sylow $19$-subgroups of $\SL$ are exactly all subgroups of $\SL$ with order $19$.
\end{lem}
\begin{proof}
Suppose that $G$ is a subgroup in $\SL$ with order $19$. Since $19$ is prime, then $G$ is a $19$-subgroup of $\SL$. Suppose that $G$ is a subgroup of some $19$-subgroup $H \subseteq \SL$. By Lagrange's Theorem, $|H|$ is a factor of $|\SL|$. \\

Since $|\SL|$ only contains one factor of $19$, then we must have ${|H| = 19}$. This implies that $G=H$. Hence, $G$ is a maximal $19$-subgroup of $\SL$ so it is a Sylow $19$-subgroup of $\SL$. This proves that all subgroups of $\SL$ with order $19$ are Sylow $19$-subgroups. \\

Since $|\SL|$ only contains one factor of $19$ then there are no Sylow $19$-subgroups of order $19^k$ for some $k \geq 2$, so this completes the proof.
\end{proof}

\begin{lem}\label{NoOrder3NonEigen}
No matrices in $\SL$ with no eigenvalues in $\F_7$ have order $1$ or $3$.
\end{lem}
\begin{proof}
The order one case is ruled out because $\lambda=1$ is an eigenvalue of the identity matrix. Suppose that some $M \in \SL$ with no eigenvalues in $\F_7$ has order $3$. Then by Lemma \ref{OrderNIffExtSame}, $M$ will have three distinct eigenvalues in $\F_{7^3}$, all with order $3$. The eigenvalues will be the roots of the equation 
$$
\lambda^3-1=0.
$$
However, the three roots are $\lambda = 1, 2$ and $4$, which is a contradiction since all the roots are all in $\F_7$.
\end{proof}

\subsection{Conjugacy Classes Preliminary Results}

\begin{lem}\label{ConjClassTrace}
Let $K$ be a conjugacy class in $\SL$. Then all matrices in $K$ have the same trace.
\end{lem}
\begin{proof}
Let $A, \ B \in K$. Then there exists a $X \in \SL$ such that $B = XAX^{-1}$. Since $\tr(NM) = \tr(MN)$ for all matrices $N,M \in \SL$, then 
$$\tr(B) = \tr((XAX^{-1}) = \tr((XA)X^{-1})= \tr(X^{-1}(XA)) = \tr((X^{-1}X)A)= \tr(A).$$
\end{proof}

\begin{lem}\label{SameConjOrder}
Let $K$ be a conjugacy class in $\SL$. Then all the matrices in $K$ have the same order.
\end{lem}
\begin{proof}
Let $M$ be a fixed matrix in $K$. Since $\SL$ is finite, then $M$ has a finite order. Consider any $A \in K$. It suffices to prove that $A, M$ have the same order. We can write $A=XMX^{-1}$ for some $X \in \SL$. For all positive integers $k$, $A^k = XM^k X^{-1}$. Then $A^k = I$ if and only if $I = XM^k X^{-1}$. This is equivalent to
$$
M^k = X^{-1} I X = I.
$$
Hence, $A$ and $M$ have the same order.
\end{proof}

\begin{lem}\label{SameConjEigen}
Let $K$ is a conjugacy class in $\SL$. Then all matrices in $K$ have the same set of eigenvalues in $\mathbb{F}_{7^3}$.
\end{lem}
\begin{proof}
Consider two matrices $A, B \in K$. Suppose that $\vv{v} \neq \vv{0}$ is an eigenvector of $A$ in $\F_{7^3}^3$ with corresponding eigenvalue $\lambda \in \F_{7^3}$. Then we have $M \vv{v} = \lambda \vv{v}$. We can write $A = XBX^{-1}$ for some $X \in \SL$. Thus, $(XBX^{-1}) \vv{v} = \lambda \vv{v}$, so ${B(X^{-1} \vv{v}) = \lambda (X^{-1} \vv{v})}$. Since $X^{-1}$ is invertible, then $X^{-1} \vv{v} \neq \vv{0}$, so $B$ also has the eigenvalue $\lambda$ (with corresponding eigenvector $X^{-1} \vv{v}$). Since $A, B$ are arbitrary elements of $K$, then this proves the result.
\end{proof}

\begin{lem}\label{NoEigenConjClass}
If a matrix $M \in \SL$ has no eigenvalues in $\F_7$, then none of the matrices in its conjugacy class have any eigenvalues in $\F_7$.
\end{lem}
\begin{proof}
This follows immediately from Lemma \ref{SameConjEigen}.
\end{proof}

\begin{lem}\label{SameConjCharacPoly}
Let $K$ be a conjugacy class in $\SL$ where none of the matrices in $K$ have eigenvalues in $\F_7$. Then all matrices $M \in K$ have the same characteristic polynomial.
\end{lem}
\begin{proof}
From Lemmas \ref{EigenVectorWhenExtended} and \ref{SameConjEigen}, all matrices $M \in K$ must have the characteristic polynomial of $p(\lambda) = (\lambda - \lambda_1)(\lambda - \lambda_2)(\lambda - \lambda_3)$ where $\lambda_1, \ \lambda_2, \ \lambda_3$ are the three distinct eigenvalues in $\F_{7^3}$ shared by all elements of $K$.
\end{proof}

\begin{lem}\label{ConjClassSizesAreSameDoubleAndQuadruple}
Let $K$ be a conjugacy class in $\SL$. Let 
\begin{align*}
    L &= \{ 2M \mid M \in K \} \subseteq \SL.\\
    T &= \{ 4M \mid M \in K \} \subseteq \SL.
\end{align*}
Then $L$ and $T$ are conjugacy classes in $\SL$ where $|K| = |L| = |T|$.
\end{lem}
\begin{proof}
The mapping $M \mapsto 2M$ defines a bijection from $K$ to $L$ so $|K| = |L|$. \\

Let $2A$ and $2B$ be arbitrary elements of $L$ where $A, B \in K$. Then, there exists a $C \in \SL$ such that $B = CAC^{-1}$. Thus, $2B = C(2A)C^{-1}$ so $2A$ and $2B$ are conjugates. Suppose that for a fixed $2A \in L$ (where $A \in K$) there exists a $D \in \SL$ such that $2A$ is a conjugate of $D$ but $D \notin L$. Then, there exists a $C \in \SL$ such that $D = C(2A)C^{-1}$. Thus, $4D = \frac{1}{2}D = CAC^{-1}$. Therefore, $4D$ and $A$ are conjugates so $4D \in K$. Hence, $8D = D \in L$, which is a contradiction. Hence, all conjugates of $2A$ must be in $L$, so $L$ is a conjugacy class. \\

Now all the elements in $T$ are of the form $2M$ where $M \in L$. The previous argument can be applied again to prove that $T$ is a conjugacy class in $\SL$ where $|L|=|T|$.
\end{proof}
\begin{lem}\label{PowerSameConjClass}
Let $M, N$ be matrices in $\SL$ that are conjugates. Then $M^k$ and $N^k$ are conjugates for all positive integers $k$.
\end{lem}
\begin{proof}
Since $M, N$ are conjugates then there exists an $X \in \SL$ such that $N = XMX^{-1}$. Then for all positive integers $k$, we have $N^k= XM^k X^{-1}$, so $M^k$ and $N^k$ are conjugates.
\end{proof}

\begin{lem}\label{OrbitStab}
Consider an $M \in \SL$. Let $G_M$ be the set of matrices in $\SL$ that commute with $M$. Then the size of the conjugacy class that contains $M$ is $|\SL|/|G_M|$.
\end{lem}
\begin{proof}
    Let $K$ be the conjugacy class that contains $M$. $K$ is a G-set under the group action of $X \cdot N = XNX^{-1}$ for all $X \in \SL$ and $N \in K$. By definition of $K$, the size of the orbit of $M$ under the group action is $|K|$. The matrices $X \in \SL $ that stabilise $M$ under the group action satisfy $XMX^{-1} = M$, which is equivalent to $XM = MX$ ($X$ commuting with $M$). The number of such matrices $X$ is $|G_M|$. By the orbit stabiliser theorem, $|\SL| = |G_M||K|$, so the result follows.
\end{proof}

\begin{lem}\label{PowerConjClassSameSize}
Let $K$ be a conjugacy class in $\SL$ where every $M \in K$ has order $n \geq 2$. For a fixed positive integer $k$ where $2 \leq k \leq n$ and $\gcd(k, n) = 1$, let 
$$
L = \{ M^k \mid M \in K \}.
$$
Then $L$ is a conjugacy class in $\SL$ where $|L| = |K|$.
\end{lem}
\begin{proof}
From Lemma \ref{PowerMatrixSetSameSizeCoprimeOrder}, $|K| = |L|$. From Lemma \ref{PowerSameConjClass}, any two elements of $|L|$ must be conjugates. Consider a fixed $A \in K$. Then $A^k \in L$. Suppose that there exists a $B \in \SL$ such that $B \notin L$ and $A^k$ and $B$ are conjugates. Then, there exists a $C \in \SL$ such that $B = CA^k C^{-1} = (CAC^{-1})^k$. Since $A \in K$, then $CAC^{-1} \in K$, so $B \in L$, which is a contradiction. Hence, $L$ contains all the conjugates of $A^k$, so $L$ is a conjugacy class.
\end{proof}

\subsection{Computations of Matrix Orders}

In the following proofs, we will use the property that
$$
(a+b)^7 = a^7+b^7
$$
for all $a, b \in \F_{7^3}$.

\begin{prop}[J.L.C]\label{43nt}
All matrices in $\mbf{[4, 3]}$ have order 19.
\end{prop}
\begin{proof} 
Since all matrices $M \in \mbf{[4, 3]}$ have no eigenvalues in $\F_7$, then by Lemma \ref{EigenVectorWhenExtended}, $M$ has three distinct eigenvalues in $\F_{7^3}$ and each eigenvalue is a root of the equation 
$$\lambda^3 - 4\lambda^2 + 3\lambda - 1 = 0.$$
By Lemma $\ref{OrderNIffExtSame}$ it suffices to prove that any root $\lambda$ of the above equation has order $19$. We will prove that $\lambda^{21} = \lambda^2$, which is equivalent to $\lambda^{19} = 1$. \\

We obtain that for all roots $\lambda$ of the cubic equation above, 
$$\lambda^3 = 4\lambda^2 - 3\lambda +1 = 4\lambda^2 + 4\lambda +1 = (2\lambda + 1)^2$$
using arithmetic in $\F_7^3$. Now we compute $\lambda^{21}$ below.

$$\begin{aligned}\lambda^{21}&=(\lambda^3)^7\\
                             &=(2\lambda+1)^{14}\\
                             &=((2\lambda+1)^7)^2\\
                             &=(2^7\lambda^7+1)^2\\
                             &=(2\lambda^7+1)^2\\
                             &=(1+2\lambda(\lambda^3)^2)^2\\
                             &=(1+2\lambda(2\lambda+1)^4)^2\\
                             &=(1+2\lambda(1+8\lambda+24\lambda^2+32\lambda^3+16\lambda^4))^2\\
                             &=(1+2\lambda(1+\lambda+3\lambda^2+4\lambda^3+2\lambda^4))^2\\
                             &=(1+2\lambda(1+\lambda+3\lambda^2+4(4\lambda^2-3\lambda+1)+2\lambda^4))^2\\
                             &=(1+2\lambda(5+3\lambda+5\lambda^2+2\lambda^4))^2\\
                             &=(1+\lambda(3-\lambda+3\lambda^2+4\lambda^4))^2\\
                             &=(1+\lambda(3-\lambda+3\lambda^2+4\lambda(4\lambda^2-3\lambda+1)))^2\\
                             &=(1+\lambda(3+3\lambda+5\lambda^2+2\lambda^3))^2\\
                             &=(1+\lambda(3+3\lambda+5\lambda^2+2(4\lambda^2-3\lambda+1)))^2\\
                             &=(1+\lambda(5+4\lambda+6\lambda^2))^2\\
                             &=(1+5\lambda+4\lambda^2-\lambda^3)^2\\
                             &=(1+5\lambda+4\lambda^2-(4\lambda^2-3\lambda+1))^2 \\
                             &=\lambda^2.\\
    \end{aligned}$$
    Hence, the order of all three eigenvalues is a factor of $19$. Since $19$ is prime and none of the eigenvalues are equal to $1$, then the three eigenvalues have order $19$.
\end{proof}

\begin{prop}[W.L]\label{31nt}
All matrices in $\mbf{[3, 1]}$ have order 19.
\end{prop}
\begin{proof}
Since all matrices $M \in \mbf{[3, 1]}$ have no eigenvalues in $\F_7$, then by Lemma \ref{EigenVectorWhenExtended}, $M$ has three distinct eigenvalues in $\F_{7^3}$ and each eigenvalue is a root of the equation 
$$\lambda^3 - 3\lambda^2 + \lambda - 1 = 0.$$

By Lemma $\ref{OrderNIffExtSame}$ it suffices to prove that any root $\lambda$ of the above equation has order $19$. Below we will compute powers of $\lambda$ to prove that $\lambda^{19} = 1$, using the relation 
$$\lambda^3 = 3\lambda^2-\lambda+1$$
which holds for all three roots $\lambda$ of the cubic equation. Below is the computation of $\lambda^4$.
 
$$\begin{aligned}\lambda^{4}&=3\lambda^{3}-\lambda^{2}+\lambda\\
&= 3(3\lambda^2-\lambda+1)-\lambda^2+\lambda \\
&=\lambda^{2}+5\lambda+3.
\end{aligned}$$
Similarly, we have:
$$
\begin{aligned}
    \lambda^{5} &=\lambda^{2}+2\lambda+1. \\
    \lambda^{6} &=5\lambda^{2}+1. \\
    \lambda^{12} &= (\lambda^6)^2. \\
                 &= (5\lambda^{2}+1)^2\\
                 &=6\lambda+6. \\
    \lambda^{18} &= \lambda^{12} \cdot \lambda^{6}.\\
                 &= (6\lambda + 6)(5\lambda^{2} + 1).\\
                 &= \lambda^{2}+4\lambda+1.\\
    \lambda^{19} &=1.\\
\end{aligned}
$$

Hence, the order of all three eigenvalues is a factor of $19$. Since $19$ is prime and none of the eigenvalues are equal to $1$, then the three eigenvalues have order $19$.
\end{proof}

\begin{prop}[L.L]\label{13nt}
All matrices in $\mbf{[1, 3]}$ have order 19.
\end{prop}

\begin{proof} 
Since all matrices $M \in \mbf{[1, 3]}$ have no eigenvalues in $\F_7$, then by Lemma \ref{EigenVectorWhenExtended}, $M$ has three distinct eigenvalues in $\F_{7^3}$ and each eigenvalue is a root of the equation 
$$\lambda^3 - \lambda^2 + 3\lambda - 1 = 0.$$

By Lemma \ref{OrderNIffExtSame} it suffices to prove that any root $\lambda$ of the above equation has order $19$. We will prove that $\lambda^{21} = \lambda^2$, which is equivalent to $\lambda^{19} = 1$. We will use the relation
$$\lambda^3 = \lambda^2 - 3\lambda +1.$$
Now we compute $\lambda^{21}$ below.
$$\begin{aligned}\lambda^{21}&=(\lambda^2-3\lambda+1)^7\\
&=\lambda^{14}-3^7\lambda^7+1\\
&=\lambda^{14}-3\lambda^7+1.\\
\end{aligned}$$
We also have
$$\begin{aligned}\lambda^6&=(\lambda^2-3\lambda+1)^2\\
&=\lambda^4+9\lambda^2+1-6\lambda^3+2\lambda^2-6\lambda \\
&= \lambda(\lambda^2-3\lambda+1) + 9\lambda^2+1-6\lambda^3+2\lambda^2-6\lambda \\
&= 2\lambda^3+\lambda^2+2\lambda+1 \\
&= 2(\lambda^2-3\lambda+1)+\lambda^2+2\lambda+1 \\
&=3(\lambda^2+\lambda+1).\\
\end{aligned}$$
This gives us
$$\begin{aligned}
\lambda^{14}&=\lambda^2(\lambda^6)^2\\
&=\lambda^2 \big(3(\lambda^2+\lambda+1) \big)^2\\
&=2\lambda^2(\lambda^4+\lambda^2+1+2\lambda^3+2\lambda^2+2\lambda)\\
&=2\lambda^6+4\lambda^5+6\lambda^4+4\lambda^3+2\lambda^2 \\
&=2(3)(\lambda^2+\lambda+1) + 4\lambda^2(\lambda^2-3\lambda+1)+6\lambda(\lambda^2-3\lambda+1)+4(\lambda^2-3\lambda+1)+2\lambda^2 \\
&=4\lambda^4-6\lambda^3-2\lambda^2+10\\
&=4\lambda(\lambda^2-3\lambda+1)-6(\lambda^2-3\lambda+1)-2\lambda^2+10\\
&=4\lambda^3-20\lambda^2+22\lambda+4\\
&=4(\lambda^2-3\lambda+1)-20\lambda^2+22\lambda+4\\
&=5\lambda^2+3\lambda+1\\
\end{aligned}$$
We now compute $-3\lambda^7$.
$$\begin{aligned}
-3\lambda^7&=-3\lambda(\lambda^6)\\
&=-3\lambda(3(\lambda^2+\lambda+1))\\
&=5\lambda^3+5\lambda^2+5\lambda\\
&=5(\lambda^2-3\lambda+1)+5\lambda^2+5\lambda\\
&=3\lambda^2-3\lambda-2\\
\end{aligned}$$
Combining the above equalities gives us
$$\begin{aligned}
\lambda^{21}&=\lambda^{14}-3\lambda^7+1\\
&=(5\lambda^2+3\lambda+1)+(3\lambda^2-3\lambda-2)+1\\
&=\lambda^2\\
\end{aligned}$$
Hence, the order of all three eigenvalues is a factor of $19$. Since $19$ is prime and none of the eigenvalues are equal to $1$, then the three eigenvalues have order $19$.
\end{proof}

\subsection{Further Results on sets of the form [i, j]}
\begin{lem}\label{OneOrderNThenAll}
Consider a set $\mbf{[i, j]} \subseteq \SL$ where no matrices in the set have any eigenvalues in $\F_7$. Then all the matrices in $\mbf{[i,j]}$ have the same order.
\end{lem}
\begin{proof}
Consider an arbitrary $M \in \mbf{[i, j]}$ with order $n$. Since $M$ has no eigenvalues in $\F_7$ then by Lemma \ref{OrderNIffExtSame}, the three distinct eigenvalues $\lambda_1, \lambda_2, \lambda_3$ of $M$ in $\F_{7^3}$ have order $n$. Consider any $N \in \mbf{[i, j]}$. Since $M$ and $N$ have the same characteristic polynomial, then $N$ also has the eigenvalues $\lambda_1, \lambda_2, \lambda_3$. By Lemma \ref{OrderNIffExtSame} again, $N$ also has order $n$.
\end{proof}

\begin{lem}\label{PowerConjClassSizeOrder19}
Consider a non-empty set ${A = \mbf{[i, j]} \subseteq \SL}$ where no matrices in the set have any eigenvalues in $\F_7$ and all matrices in the set have order $19$. For a fixed positive integer $k$ where $2 \leq k \leq 18$, let 
$$B = \{ M^k \mid M \in A \} \subseteq \SL.$$
Then $B$ is of the form $\mbf{[v,w]}$ for some $v$ and $w$, and none of the matrices in $B$ have any eigenvalues in $\F_7$. Also $|A| = |B|$ and all the matrices in $B$ have order $19$.
\end{lem}
\begin{proof}
Since $\gcd(k, 19) = 1$ and $k \neq 1$, then there exists a positive integer $t$ such that ${2 \leq t \leq 18}$ and ${\ModCmp{kt}{1}{19}}$. Thus, $kt = 19s+1$ for some non-negative integer $s$. \\

Consider any $M \in A$. Since $M$ has order $19$, then by Lemma \ref{orderPowerCoprimeFixed}, $M^k$ also has order $19$. By Lemma \ref{EigenCoprime}, none of the matrices in $B$ have any eigenvalues in $\F_7$. Since $M$ has no eigenvalues in $\F_7$, then by Lemma \ref{OrderNIffExtSame}, the characteristic polynomial for $M$ must be ${P(\lambda) = (\lambda - \lambda_1)(\lambda - \lambda_2)(\lambda - \lambda_3)}$ for distinct eigenvalues ${\lambda_1, \lambda_2, \lambda_3 \in \F_{7^3}}$ where the eigenvalues all have order $19$.\\

Note that 
$$({\lambda_1}^k)^t = {\lambda_1}^{19s+1} = \lambda_1.$$ 
Similarly, $({\lambda_2}^k)^t = \lambda_2$ and $({\lambda_3}^k)^t = \lambda_3$. Thus, since $\lambda_1, \ \lambda_2$ and $\lambda_3$ are distinct, then $\lambda_1^k, \ \lambda_2^k$ and $\lambda_3^k$ are distinct. Then by Lemma \ref{EigenvalueScale}, $\lambda_1^k, \ \lambda_2^k$ and $\lambda_3^k$ are the three distinct eigenvalues of $M^k$ where these eigenvalues are in $\F_{7^3}$ but not in $\F_7$ (because $M^k$ has no eigenvalues in $\F_7$). \\

Hence, the characteristic polynomial of $M^k$ must be $Q(\lambda) = (\lambda - \lambda_1^k)(\lambda - \lambda_2^k)(\lambda - \lambda_3^k)$. Note that all matrices in $A$ will have the same characteristic polynomial $P$ and since $Q$ is uniquely determined by $P$, then all matrices in $B$ will also have the same characteristic polynomial $Q$, which has no roots in $\F_7$. Therefore, $B \subseteq C$ where $C = \mbf{[v, w]}$ for some $v$ and $w$ and also no matrices in $C$ have eigenvalues in $\F_7$. \\

By Lemma \ref{PowerMatrixSetSameSizeCoprimeOrder}, $|A| = |B|$. Since all matrices in $B \subseteq C$ have order $19$, then by Lemma \ref{OneOrderNThenAll}, all matrices in $C$ have order $19$. Let
$$
D = \{ M^t \mid M \in C \} \subseteq \SL.
$$
Since an arbitrary element of $B$ is of the form $M^k$ for some $M \in A$ and since $B \subseteq C$, then $(M^k)^t \in D$. We can compute 
$$(M^k)^t = M^{19s+1} = M^{19s} \cdot M = M.$$ 
Thus, $M \in D$, so $A \subseteq D$. \\

Using a similar argument as earlier, we also obtain that $D \subseteq \mbf{[x, y]}$ for some $x$ and $y$ and $|C| = |D|$. Then, we have 
$$
\mbf{[i, j]} = A \subseteq D \subseteq \mbf{[x, y]}.
$$
This forces us to have ${\mbf{[i, j]} = \mbf{[x, y]}}$ and so we must have $A = D$, which gives us ${|C| = |D| = |A| = |B|}$. Since $B \subseteq C$, then we must have $B = C$. Hence, $B = \mbf{[v, w]}$.
\end{proof} 

\begin{lem}\label{AllPairsHaveTheSameSize}
The sets $\mbf{[0, 1]}$, $\mbf{[0, 2]}$, $\mbf{[0, 4]}$, $\mbf{[1, 0]}$, $\mbf{[1, 3]}$, $\mbf{[1, 5]}$, $\mbf{[2, 0]}$, $\mbf{[2, 5]}$, $\mbf{[2, 6]}$, $\mbf{[3, 1]}$, $\mbf{[3, 4]}$, $\mbf{[4, 0]}$, $\mbf{[4, 3]}$, $\mbf{[4, 6]}$, $\mbf{[5, 1]}$, $\mbf{[5, 2]}$, $\mbf{[6, 2]}$, $\mbf{[6, 4]}$ all have the same size.
\end{lem}
\begin{proof}
We will first apply Lemma \ref{M2MBijection} repeatedly, as shown in the following equalities. 
\begin{align*} 
|\mbf{[0,1]}| &= |\mbf{[0, 4]}| = |\mbf{[0,2]}|.\\
|\mbf{[1,0]}| &= |\mbf{[2,0]}| = |\mbf{[4,0]}|.\\
|\mbf{[1,3]}| &= |\mbf{[2,5]}| = |\mbf{[4,6]}|.\\
|\mbf{[1,5]}| &= |\mbf{[2,6]}| = |\mbf{[4,3]}|.\\
|\mbf{[3,1]}| &= |\mbf{[6,4]}| = |\mbf{[5,2]}|.\\
|\mbf{[3,4]}| &= |\mbf{[6,2]}| = |\mbf{[5,1]}|.\\
\end{align*}
Let 
$$
M = \MatrThree{0}{2}{-1}{0}{0}{2}{2}{0}{0}.
$$
We proved in the table in Section 2 that $M$ has order $19$. Also from the table we have the following:

\begin{align*}
M &\in  \mbf{[0,2]}.\\
M^2 &\in \mbf{[0,2]}.\\
M^4 &\in \mbf{[0,2]}.\\
M^5 &\in \mbf{[0,2]}.\\
M^8 &\in \mbf{[0,2]}.\\
M^9 &\in \mbf{[3,1]}.\\
\end{align*}
By Lemma \ref{OneOrderNThenAll} and \ref{18PairsRepresentAllNonEigenMatrices}, all matrices in $\mbf{[0,2]}$ must have order $19$ and none have any eigenvalues in $\F_7$. Then using Lemma \ref{PowerConjClassSizeOrder19} we get that 
$$
|\mbf{[0,2]}| = |\mbf{[3, 4]}| = |\mbf{[1,3]}| = |\mbf{[4, 3]}| = |\mbf{[2,0]}| = |\mbf{[3,1]}|.
$$
Combining this with the earlier equalities gives us that all $18$ sets have the same size.
\end{proof}

\subsection{Further Results on Conjugacy Classes}

\begin{lem}\label{PowerOfMInConjClass}
Consider an $M \in \SL$ with order $19$. Consider a collection of conjugacy classes ${K_1, K_2, \ldots, K_m}$ in $\SL$ where all the elements in these conjugacy classes have order $19$. Suppose that $M \in K_1$. Then for each integer $n$ where ${1 \leq n \leq m}$, there exists a positive integer $c_n$ where ${1 \leq c_n \leq 18}$ such that $M^{c_n} \in K_n$. 
\end{lem}
\begin{proof}
Consider some value of $n$ where ${1 \leq n \leq m}$. Let $P$ be the cyclic group of order $19$ generated by $M$. Choose an arbitrary $N \in K_n$ (which has order $19$ by definition). Let $Q$ be the cyclic group of order $19$ generated by $N$. By Lemma \ref{AllOrder19SubgroupsAreSylow}, $P$ and $Q$ are Sylow $19$-subgroups of $\SL$. By Sylow's second theorem, $P$ and $Q$ are conjugate to each other, so there exists an $A \in \SL$ such that $A^{-1}QA = P$. Hence, there exists a positive integer $t$ where $1 \leq t \leq 18$ such that $A^{-1}NA = M^t$. Since $K_m$ is a conjugacy class, then 
$$
M^t = A^{-1}NA \in K_n.
$$
\end{proof}


\begin{lem}\label{ConjClassWithinAPair}
Let $K$ be a conjugacy class in $\SL$ where none of the matrices in $K$ have eigenvalues in $\F_7$. Then, $K \subseteq \mbf{[i, j]}$ for some values of $i$ and $j$.
\end{lem}
\begin{proof}
This immediately follows from Lemma \ref{SameConjCharacPoly}.
\end{proof}

\begin{lem}\label{DistinctConjClassesSameSizeFoundInThreePairs}
There exist distinct conjugacy classes in $\SL$ contained in each of $\mbf{[0, 1]}$, $\mbf{[0, 2]}$ and $\mbf{[0, 4]}$, all of which have size $2^5 \times 3^2 \times 7^3$.
\end{lem}
\begin{proof}
Let
\begin{align*}
    M_1 &= \MatrThree{0}{1}{3}{0}{0}{1}{1}{0}{0}.\\
    M_2 &= 2M_1.\\
    M_3 &= 4M_1.\\
\end{align*}
We saw in Section $1$ that $M_1 \in \mbf{[0, 4]}$. By Lemma \ref{M2MBijection} we get that $M_2 \in \mbf{[0, 2]}$. By applying Lemma \ref{M2MBijection} on $\mbf{[0, 2]}$ we obtain $M_3 \in \mbf{[0, 1]}$. \\

Let $K_1, K_2, K_3$ be the conjugacy classes in $\SL$ which contain $M_1, M_2, M_3$ respectively. By Lemma \ref{ConjClassWithinAPair}, we have 
\begin{align*}
K_1 &\subseteq \mbf{[0, 4]}.\\
K_2 &\subseteq \mbf{[0, 2]}.\\
K_3 &\subseteq \mbf{[0, 1]}.\\
\end{align*}
In Section 3 we showed that ${|K_2| = 2^5 \times 3^2 \times 7^3}$. By Lemma  \ref{ConjClassSizesAreSameDoubleAndQuadruple}, $K_1$, $K_2$ and $K_3$ all have the same size of ${2^5 \times 3^2 \times 7^3}$.
\end{proof}

\begin{lem}\label{FiveMoreConjClassesFound}
There exist distinct conjugacy classes in $\SL$ contained within each of $\mbf{[1, 3]}$, $\mbf{[2, 0]}$, $\mbf{[3, 1]}$, $\mbf{[3, 4]}$, $\mbf{[4, 3]}$, all of which have size $2^5 \times 3^2 \times 7^3$.
\end{lem}
\begin{proof}
For $i=1, 2, 3$ define $M_i$ and $K_i$ as in Lemma \ref{DistinctConjClassesSameSizeFoundInThreePairs}. Note that the matrix 
$$M_2=2M_1= \MatrThree{0}{2}{-1}{0}{0}{2}{2}{0}{0}.$$
was shown to have order $19$ from the table in Section $2$. Let
\begin{align*}
    M_4 &= {M_2}^6.\\
    M_5 &= {M_2}^8.\\
    M_6 &= {M_2}^{15}.\\
    M_7 &= {M_2}^2.\\
    M_8 &= {M_2}^5.\\
\end{align*}

 From that same table, we have 
 \begin{align*}
     M_4 &\in \mbf{[1, 3]}.\\
     M_5 &\in \mbf{[2,0]}.\\
     M_6 &\in \mbf{[3, 1]}.\\
     M_7 &\in \mbf{[3, 4]}.\\
     M_8 &\in \mbf{[4, 3]}.\\
 \end{align*}
 For $t = 4, 5, \ldots, 8$, let $K_t$ be the conjugacy class in $\SL$ that contains $M_t$. By Lemma \ref{ConjClassWithinAPair}, we have
 \begin{align*}
     K_4 &\subseteq \mbf{[1, 3]}.\\
     K_5 &\subseteq \mbf{[2, 0]}.\\
     K_6 &\subseteq \mbf{[3, 1]}.\\
     K_7 &\subseteq \mbf{[3, 4]}.\\
     K_8 &\subseteq \mbf{[4, 3]}.\\
 \end{align*}
From applying Lemma \ref{PowerConjClassSameSize} to $M_2$ (which has order $19$), we obtain that $K_4, K_5, \ldots, K_8$ have the same size as $K_2$, which we know to be ${2^5 \times 3^2 \times 7^3}$ from Lemma \ref{DistinctConjClassesSameSizeFoundInThreePairs}.
\end{proof}

\begin{lem}\label{RemainingConjClassSizesFound}
There exist distinct conjugacy classes in $\SL$ contained within each of $\mbf{[2, 5]}$, $\mbf{[4, 6]}$, $\mbf{[4, 0]}$, $\mbf{[1, 0]}$, $\mbf{[6, 4]}$, $\mbf{[5, 2]}$, $\mbf{[6, 2]}$, $\mbf{[5, 1]}$, $\mbf{[1, 5]}$, $\mbf{[2, 6]}$, all of which have size $2^5 \times 3^2 \times 7^3$.
\end{lem}
\begin{proof}
For $i=1, 2, \ldots, 8$, define $M_i$ and $K_i$ as in Lemma \ref{FiveMoreConjClassesFound}. Then, we will define
\begin{align*}
M_9 &= 2M_4.\\
M_{10} &= 4M_4.\\
M_{11} &= 2M_5.\\
M_{12} &= 4M_5.\\
M_{13} &= 2M_6.\\
M_{14} &= 4M_6.\\
M_{15} &= 2M_7.\\
M_{16} &= 4M_7.\\
M_{17} &= 2M_8.\\
M_{18} &= 4M_8.\\
\end{align*}
For $9 \leq i \leq 18$, let $K_i$ be the conjugacy class in $\SL$ that contains $M_i$. From applying Lemma \ref{M2MBijection} several times, we have 
\begin{align*}
M_9 &\in \mbf{[2, 5]}.\\
M_{10} &\in \mbf{[4, 6]}.\\
M_{11} &\in \mbf{[4, 0]}.\\
M_{12} &\in \mbf{[1, 0]}.\\
M_{13} &\in \mbf{[6, 4]}.\\
M_{14} &\in \mbf{[5, 2]}.\\
M_{15} &\in \mbf{[6, 2]}.\\
M_{16} &\in \mbf{[5, 1]}.\\
M_{17} &\in \mbf{[1, 5]}.\\
M_{18} &\in \mbf{[2, 6]}.\\
\end{align*}
By Lemma \ref{ConjClassWithinAPair},  we have 
\begin{align*}
K_9 &\subseteq \mbf{[2, 5]}.\\
K_{10} &\subseteq \mbf{[4, 6]}.\\
K_{11} &\subseteq \mbf{[4, 0]}.\\
K_{12} &\subseteq \mbf{[1, 0]}.\\
K_{13} &\subseteq \mbf{[6, 4]}.\\
K_{14} &\subseteq \mbf{[5, 2]}.\\
K_{15} &\subseteq \mbf{[6, 2]}.\\
K_{16} &\subseteq \mbf{[5, 1]}.\\
K_{17} &\subseteq \mbf{[1, 5]}.\\
K_{18} &\subseteq \mbf{[2, 6]}.\\
\end{align*}
It was proven in Lemma \ref{FiveMoreConjClassesFound} that $K_4, K_5, \ldots, K_8$ each have size ${2^5 \times 3^2 \times 7^3}$. From applying Lemma \ref{ConjClassSizesAreSameDoubleAndQuadruple} to $M_4, M_5, \ldots, M_8$, we obtain that $K_9, K_{10}, \ldots, K_{18}$ also have size ${2^5 \times 3^2 \times 7^3}$.
\end{proof}
\begin{lem}\label{All18ConjClassSizesFound}
There exist at least $18$ distinct conjugacy classes in $\SL$ of size $2^5 \times 3^2 \times 7^3$ where each is contained in one of $\mbf{[0, 1]}$, $\mbf{[0, 2]}$, $\mbf{[0, 4]}$, $\mbf{[1, 0]}$, $\mbf{[1, 3]}$, $\mbf{[1, 5]}$, $\mbf{[2, 0]}$, $\mbf{[2, 5]}$, $\mbf{[2, 6]}$, $\mbf{[3, 1]}$, $\mbf{[3, 4]}$, $\mbf{[4, 0]}$, $\mbf{[4, 3]}$, $\mbf{[4, 6]}$, $\mbf{[5, 1]}$, $\mbf{[5, 2]}$, $\mbf{[6, 2]}$ and $\mbf{[6, 4]}$.
\end{lem}

\begin{proof}
This follows from Lemmas \ref{DistinctConjClassesSameSizeFoundInThreePairs}, \ref{FiveMoreConjClassesFound} and \ref{RemainingConjClassSizesFound}.
\end{proof}

\begin{lem}\label{CharacteriseAllOrder19And57Matrices}
All the matrices in $\mbf{[0,2]}$, $\mbf{[1,3]}$, $\mbf{[2,0]}$,  $\mbf{[3,1]}$, $\mbf{[3,4]}$ and $\mbf{[4,3]}$ have order $19$. Also, all the matrices in $\mbf{[0,4]}$, $\mbf{[0,1]}$, $\mbf{[2,5]}$, $\mbf{[4,6]}$, $\mbf{[4,0]}$, $\mbf{[1,0]}$, $\mbf{[6,4]}$, $\mbf{[5,2]}$, $\mbf{[6,2]}$, $\mbf{[5,1]}$, $\mbf{[1,5]}$, $\mbf{[2,6]}$ have order $57$.
\end{lem}
\begin{proof}
For $i=1, 2, 3, \ldots, 18$ define $M_i$ and $K_i$ as in Lemma \ref{RemainingConjClassSizesFound}. Recall that $M_2 \in \mbf{[0, 2]}$ has order $19$. For $j= 4, 5, \ldots, 8$ recall that $M_j$ was defined to be of the form $M_2^k$ for certain values of $k$ where $2 \leq k \leq 18$. Hence, each such $M_j$ has order $19$ because the values of $k$ are coprime to $19$. \\

Next, recall that for $t=9, 10, \ldots, 18$ the matrix $M_t$ was defined to be of the form $2M_j$ or $4M_j$ for certain values of $j$ where $4 \leq j \leq 8$. Also, we have
$$M_1 = \frac{M_2}{2}=4M_2$$
and $M_3=2M_2$. By Lemma \ref{OrderOfDoubleAndQuadrupleIs57} we obtain that the matrices $M_1, M_3$ and $M_t$ for $t=9, 10, \ldots, 18$ all have order $57$. Thus, we have found representative matrices with order $19$ in each of $\mbf{[0,2]}$, $\mbf{[1,3]}$, $\mbf{[2,0]}$,  $\mbf{[3,1]}$, $\mbf{[3,4]}$ and $\mbf{[4,3]}$ and representative matrices with order $57$ in each of $\mbf{[0,4]}$, $\mbf{[0,1]}$, $\mbf{[2,5]}$, $\mbf{[4,6]}$, $\mbf{[4,0]}$, $\mbf{[1,0]}$, $\mbf{[6,4]}$, $\mbf{[5,2]}$, $\mbf{[6,2]}$, $\mbf{[5,1]}$, $\mbf{[1,5]}$, $\mbf{[2,6]}$. \\

The result then follows from Lemma \ref{OneOrderNThenAll}.


\end{proof}

\begin{lem}\label{AllMatricesNonEigenOrder19Or57}
All matrices in $\SL$ without any eigenvalues in $\F_7$ have order $19$ or $57$.
\end{lem}
\begin{proof}
There are two ways to arrive at this result. One way is to apply Lemmas \ref{18PairsRepresentAllNonEigenMatrices} and \ref{CharacteriseAllOrder19And57Matrices} together. The second way is to apply Lemmas \ref{OrderNIffExtSame} and \ref{NoOrder3NonEigen} together (since the positive integers that divide $57$ are $1, 3, 19, 57$).
\end{proof}

\begin{lem}\label{AllOrder19MatricesNoEigen}
All matrices in $\SL$ with order $19$ have no eigenvalues in $\F_7$.
\end{lem}
\begin{proof}
Define $M_2$ and $K_2$ as in Lemma \ref{DistinctConjClassesSameSizeFoundInThreePairs}. We know from earlier that $M_2$ has no eigenvalues in $\F_7$ and has order $19$. Also, $K_2$ is the conjugacy class that contains $M_2$. Let $N$ be an arbitrary matrix in $\SL$ with order $19$ and let $K$ be the conjugacy class that contains it.\\

By Lemma \ref{PowerOfMInConjClass}, there exists a positive integer $c$ where $1 \leq c \leq 18$ such that $M_2^c \in K$. Since $19$ is prime, then by Lemma \ref{PrimeNoEigen}, $M_2^c$ has no eigenvalues in $\F_7$. Then by Lemma \ref{NoEigenConjClass}, all matrices in $K$ have no eigenvalues in $\F_7$, so $N$ has no eigenvalues in $\F_7$. 
\end{proof}

\begin{lem}\label{AllOrder57MatricesNoEigen}
All matrices in $\SL$ with order $57$ have no eigenvalues in $\F_7$.
\end{lem}
\begin{proof}
Suppose the contrary, that a matrix $M \in \SL$ has order $57$ and an eigenvalue $\lambda \in \F_7$. Then $M^3$ has order $19$ and also has eigenvalue $\lambda^3 \in \F_7$ by Lemma \ref{EigenvalueScale}. This contradicts Lemma \ref{AllOrder19MatricesNoEigen}.
\end{proof}

\begin{lem}\label{NoOtherOrder19Or57Matrices}
The matrices given in Lemma \ref{CharacteriseAllOrder19And57Matrices} describe all the matrices in $\SL$ that have order $19$ or $57$.
\end{lem}
\begin{proof}
This follows from Lemmas \ref{AllOrder19MatricesNoEigen} and \ref{AllOrder57MatricesNoEigen} because the matrices given in Lemma \ref{CharacteriseAllOrder19And57Matrices} cover all possible matrices in $\SL$ without any eigenvalues in $\F_7$.
\end{proof}

\begin{lem}\label{Matrices18And54TimesSubgroups}
Let $A$ be the set of all matrices in $\SL$ with order $19$. Let $B$ be the set of matrices in $\SL$ with order $19$ or $57$. Using the notation given earlier in the stating of Sylow's theorems, let $n_{19}$ be the number of Sylow $19$-subgroups of $\SL$ (which we know to be precisely the subgroups of order $19$ from Lemma \ref{AllOrder19SubgroupsAreSylow}). Then ${|A| = 18n_{19}}$ and ${|B| = 54n_{19}}$.

\end{lem}

\begin{proof}
All the Sylow $19$-subgroups are cyclic since $19$ is prime. The non-identity elements of all of these subgroups have order $19$ and so are in $A$. Also, the subgroups generated by each element of $A$ are cyclic subgroups of order $19$ and thus are all Sylow $19$-subgroups. Hence, using Lemma \ref{DifferentSubgroupsDisjoint}, we obtain that $A$ can be partitioned into subsets of size $18$ where each subset contains the $18$ non-identity elements of each Sylow $19$-subgroup. Hence, ${|A| = 18|n_{19}|}$. \\

Let $m = |\mbf{[0,2]}|$. By Lemmas \ref{18PairsRepresentAllNonEigenMatrices}, \ref{AllPairsHaveTheSameSize}, \ref{CharacteriseAllOrder19And57Matrices} and \ref{NoOtherOrder19Or57Matrices}, it follows that ${|A| = 6m}$ and ${|B| = 18m}$. Hence,
$$|B| = 3|A| = 54n_{19}.$$ 

\end{proof}

\begin{lem}\label{PossibleSubgroupSize}
 Define $n_{19}$ as in Lemma \ref{Matrices18And54TimesSubgroups}. Then,
 $$n_{19} \in \{1, \  2^5 \cdot 3, \ 7^3, \ 2^4 \cdot 3^2 \cdot 7, \ 2^3 \cdot 3^3 \cdot 7^2, \ 2^5 \cdot 3 \cdot 7^3       \}.$$
\end{lem} 
\begin{proof}
The order of $\SL$ can be written as $19m$ where $m = 2^5 \times 3^3 \times 7^3$. By Sylow's third theorem with the prime $p = 19$, we obtain the following conditions on $n_{19}$.

\begin{itemize}
    \item $n_{19} \equiv 1 \; (\bmod\; 19)$
    \item $n_{19}$ is a factor of $m$.
\end{itemize}
By testing the out all factors of $m$ that are $1 \; (\bmod\; 19)$, we obtain the possible values of $n_{19}$ listed above.
\end{proof}

\begin{lem}\label{The18PairsThemselvesAreConjClasses}
There exists exactly $18$ distinct conjugacy classes in $\SL$ where none of the matrices in the conjugacy classes have eigenvalues in $\F_7$. Each conjugacy class has size $2^5 \times 3^2 \times 7^3$. These are the sets $\mbf{[0, 1]}$, $\mbf{[0, 2]}$, $\mbf{[0, 4]}$, $\mbf{[1, 0]}$, $\mbf{[1, 3]}$, $\mbf{[1, 5]}$, $\mbf{[2, 0]}$, $\mbf{[2, 5]}$, $\mbf{[2, 6]}$, $\mbf{[3, 1]}$, $\mbf{[3, 4]}$, $\mbf{[4, 0]}$, $\mbf{[4, 3]}$, $\mbf{[4, 6]}$, $\mbf{[5, 1]}$, $\mbf{[5, 2]}$, $\mbf{[6, 2]}$ and $\mbf{[6, 4]}$.
\end{lem}
\begin{proof}
Define $M_1, \ M_2, \ldots, \ M_{18}$ and $K_1, \ K_2, \ldots, \ K_{18}$ as in Lemmas \ref{DistinctConjClassesSameSizeFoundInThreePairs}, \ref{FiveMoreConjClassesFound} and \ref{RemainingConjClassSizesFound}. Let $m$ be the total number of matrices in $\SL$ without any eigenvalues in $\F_7$. Also define $A$, $B$ and $n_{19}$ as in Lemma \ref{Matrices18And54TimesSubgroups}. \\

From Lemma \ref{NoOtherOrder19Or57Matrices}, we obtain that $|B| = m$. In  Lemma \ref{All18ConjClassSizesFound} we proved that there exists a conjugacy class contained within each of the $18$ sets listed above, where each of the conjugacy classes has size ${2^5 \times 3^2 \times 7^3}$. \\

It follows that 
$$|B| = m \geq 2^5 \times 3^2 \times 7^3 \times 18 = 2^6 \times 3^4 \times 7^3.$$ 
From Lemma \ref{Matrices18And54TimesSubgroups}, we have $|B| = 54n_{19}$ so 
$$n_{19} \geq \frac{2^6 \times 3^4 \times 7^3}{54} = 2^5 \times 3 \times 7^3.$$ 
Then from Lemma \ref{PossibleSubgroupSize}, we must have $n_{19}=2^5 \times 3 \times 7^3$. Thus, equality must occur in the previous inequalities so we must have $m = |B| = 2^6 \times 3^4 \times 7^3$. This also means that the $18$ sets listed above themselves must be conjugacy classes. \\

Since these $18$ conjugacy classes partition the set of matrices in $\SL$ with no eigenvalues in $\F_7$, then these form all the conjugacy classes in $\SL$ where none of the matrices in the conjugacy classes have eigenvalues in $\F_7$.
\end{proof}

\begin{lem}\label{TotalNumberOfNonEigenMatricesFound}
There are exactly $2^6 \times 3^4 \times 7^3$ matrices in $\SL$ without any eigenvalues in $\F_7$.
\end{lem}
\begin{proof}
This result follows from the proof of Lemma \ref{The18PairsThemselvesAreConjClasses}.
\end{proof}

\newpage
\subsection{Relationships between classes}

Here is a table of the 20 powers of a matrix with trace 1
\ \\

\begin{tabular}{c|c|c|c}
power&matrix&trace&class\\\hline
1&$\begin{pmatrix} 
0&1&-3\\
0&1&-2\\
1&0&0\\
\end{pmatrix}$&1&$\mbf{[1,3]}$\\\hline
2&$\begin{pmatrix} 
-3&1&-2\\
-2&1&-2\\
0&1&-3\\
\end{pmatrix}$ &2&$\mbf{[2,0]}$\\\hline
3& $\begin{pmatrix} 
-2&-2&0\\
-2&-1&-3\\
-3&1&-2\\
\end{pmatrix}$&2&$\mbf{[2,0]}$\\\hline
4&$\begin{pmatrix} 
0&3&3\\
-3&-3&1\\
-2&-2&0\\
\end{pmatrix}$&4&$\mbf{[4,3]}$\\\hline
5&$\begin{pmatrix} 
3&3&1\\
1&1&1\\
0&3&3\\
\end{pmatrix}$&0&$\mbf{[0,2]}$\\\hline
6& $\begin{pmatrix} 
1&-1&-1\\
1&2&2\\
3&3&1\\
\end{pmatrix}$&4&$\mbf{[4,3]}$\\\hline
7&$\begin{pmatrix} 
-1&0&-1\\
2&3&0\\
1&-1&-1\\
\end{pmatrix}$ &1&$\mbf{[1,3]}$\\\hline
8& $\begin{pmatrix} 
-1&-1&3\\
0&-2&2\\
1&0&0\\
\end{pmatrix}$&3&$\mbf{[3,1]}$\\\hline
9& $\begin{pmatrix} 
3&-2&-2\\
2&-2&-3\\
-1&-1&3\\
\end{pmatrix}$&4&$\mbf{[4,3]}$\\\hline
10& $\begin{pmatrix} 
-2&1&2\\
-3&0&-2\\
3&-2&-2\\
\end{pmatrix}$&3&$\mbf{[3,4]}$\\\hline
\end{tabular}
\quad
\begin{tabular}{c|c|c|c}
power&matrix&trace&class\\\hline
11&$\begin{pmatrix} 
2&-1&-3\\
-2&-3&2\\
-2&1&2\\
\end{pmatrix}$ &1&$\mbf{[1,3]}$\\\hline
12&$\begin{pmatrix} 
-3&1&3\\
2&2&-2\\
2&-1&-3\\
\end{pmatrix}$ &3&$\mbf{[3,1]}$\\\hline
13&$\begin{pmatrix} 
3&-2&0\\
-2&-3&-3\\
-3&1&3\\
\end{pmatrix}$ &3&$\mbf{[3,4]}$\\\hline
14&$\begin{pmatrix} 
0&1&2\\
-3&2&-2\\
3&-2&0\\
\end{pmatrix}$ &2&$\mbf{[2,0]}$\\\hline
15&$\begin{pmatrix} 
2&1&-2\\
-2&-1&-2\\
0&1&2\\
\end{pmatrix}$ &3&$\mbf{[3,4]}$\\\hline
16&$\begin{pmatrix} 
-2&3&-1\\
-2&-3&1\\
2&1&-2\\
\end{pmatrix}$ &0&$\mbf{[0,2]}$\\\hline
17&$\begin{pmatrix} 
-1&1&0\\
1&2&-2\\
-2&3&-1\\
\end{pmatrix}$ &0&$\mbf{[0,2]}$\\\hline
18&$\begin{pmatrix} 
0&0&1\\
-2&3&0\\
-1&1&0\\
\end{pmatrix}$ &3&$\mbf{[3,1]}$\\\hline 
19&$\begin{pmatrix} 
1&0&0\\
0&1&0\\
0&0&1\\
\end{pmatrix}$&3&$\mbf{[3,3]}$\\\hline
20&$\begin{pmatrix} 
0&1&-3\\
0&1&-2\\
1&0&0\\
\end{pmatrix}$&1&$\mbf{[1,3]}$\\\hline
\end{tabular} \\

We have bijections
$$\mbf{[1,3]}\to\mbf{[0,2]}$$
$$M\mapsto M^5$$
with inverse $$M\mapsto M^{4}.$$
Also $$M\mapsto M^{16}$$
with inverse $$M\mapsto M^{6}.$$
Also $$M\mapsto M^{17}$$
with inverse $$M\mapsto M^{9}.$$

This can be explained as the three bijections from (1,3) to itself composed with a single bijection to another class. The three bijections are identity, 7th power and $7^2=11$ th power.

\subsection{Subgroups section}
\begin{lem}\label{MaximalSubgroup}
Let $H$ be a proper subgroup of $G = \SL$. Let $A = \{A_1, A_2, \ldots A_n\}$ be a generating set for $G$. Suppose that for any $M \in G/H$, there exists matrices $B_1, B_2, \ldots B_n, C_1, C_2, \ldots C_n \in H$ such that $B_m M C_m = A_m$ for $m = 1, 2, \ldots n$. Then $H$ is a maximal subgroup of $G$.
\end{lem}
\begin{proof}
Suppose the contrary, where there exists a proper subgroup $K$ of $G$ such that $H$ is a proper subgroup of $K$. Then, consider a matrix $M \in K/H$. Since $K$ is a proper subgroup of $G$, then there exists a generator $X$ in $A$ that is not in $K$. However, there exists matrices $B$ and $C$ in $H$ such that $BMC = X$. Since $B, M$ and $C$ are in $K$, then $X$ is in $K$, contradiction. Thus, $H$ is a maximal subgroup of $G$. 
\end{proof}

\begin{thm}The subgroup $H$ of matrices of the form
$$\begin{pmatrix}
a&b&c\\
0&e&f\\
0&h&i\\
\end{pmatrix}$$
is maximal in $\SL$. 
\end{thm}
\begin{proof}

The generators of $\SL$ are as follows:

$$X = \begin{pmatrix}
1&0&1\\
0&-1&-1\\
0&1&0\\
\end{pmatrix}$$

$$Y = \begin{pmatrix}
0&1&0\\
0&0&1\\
1&0&0\\
\end{pmatrix}$$

$$Z = \begin{pmatrix}
0&1&0\\
1&0&0\\
-1&-1&-1\\
\end{pmatrix}$$\\

To prove that $H$ is a maximal subgroup, by \ref{EigenVectorWhenExtended}, it suffices to prove that given any arbitrary matrix $A \in \SL$ that is not in $H$, the generators $X, Y, Z$ for $\SL$ can be obtained by multiplying $A$ (on the left and the right) only by elements of $H$. Since $X$ itself is an element of $H$, we only need to prove that generators $Y$ and $Z$ are obtainable. \\

\subsubsection{Generators Proof for $Y$}
Let 

$$
A = \begin{pmatrix}
a&b&c\\
d&e&f\\
g&h&i\\
\end{pmatrix}
$$\\

Since $A$ is not in $H$, then either $d$ or $g$ must be non-zero. In the case of $g = 0$, then $d$ must be non-zero. We can multiply a matrix on the left to add row $2$ of $A$ to row $3$ to obtain the matrix $B$ as shown below.

$$B = \begin{pmatrix}
1&0&0\\
0&1&0\\
0&1&1\\
\end{pmatrix}
\begin{pmatrix}
a&b&c\\
d&e&f\\
0&h&i\\
\end{pmatrix} = 
\begin{pmatrix}
a&b&c\\
d&e&f\\
d&e+h&f+i\\
\end{pmatrix}$$\\

Otherwise, if $d = 0$ and $g \neq 0$, then we can multiply a matrix on the left to add row $3$ of $A$ to row $2$ to obtain the matrix $B$ as shown below.

$$ B = \begin{pmatrix}
1&0&0\\
0&1&1\\
0&0&1\\
\end{pmatrix}
\begin{pmatrix}
a&b&c\\
0&e&f\\
g&h&i\\
\end{pmatrix} = 
\begin{pmatrix}
a&b&c\\
g&e+h&f+i\\
g&h&i\\
\end{pmatrix}$$\\

Finally, if both $d$ and $g$ in $A$ were non-zero, then we can set $B =A$. In all cases we can obtain a matrix $B$ of the form 

$$ B = \begin{pmatrix}
a&b&c\\
d&e&f\\
g&h&i\\
\end{pmatrix}$$\\

where $d$ and $g$ are non-zero. Now, we can multiply a matrix on the left of $B$ to add multiples of row $3$ to row $2$ and row $1$. We can multiply by another matrix with determinant $1$ to scale row $3$. This procedure can yield a matrix $C$ with its top-left and left-middle entries being zero and its bottom-left entry being $1$, as shown below.

$$ C = 
\begin{pmatrix}
1&0&0\\
0&\frac{1}{l}&0\\
0&0&l\\
\end{pmatrix}
\begin{pmatrix}
1&0&j\\
0&1&k\\
0&0&1\\
\end{pmatrix}
\begin{pmatrix}
a&b&c\\
d&e&f\\
g&h&i\\
\end{pmatrix} = 
\begin{pmatrix}
0&m&n\\
0&p&q\\
1&s&t\\
\end{pmatrix}$$\\

Write $C$ in the form 

$$ C = \begin{pmatrix}
0&b&c\\
0&e&f\\
1&h&i\\
\end{pmatrix}$$\\

Now we can multiply a matrix on the right of $C$ to add multiples of column $1$ to columns $2$ and $3$ to yield matrix $D$ with its bottom-middle and bottom-right entries zero.

$$ D = \begin{pmatrix}
0&b&c\\
0&e&f\\
1&h&i\\
\end{pmatrix}
\begin{pmatrix}
1&-h&-i\\
0&1&0\\
0&0&1\\
\end{pmatrix} = 
\begin{pmatrix}
0&b&c\\
0&e&f\\
1&0&0\\
\end{pmatrix}$$\\

Note that $D \in \SL$ because we performed multiplications only with matrices with determinant $1$. Hence, $\mathrm{det}(D) = bf-ce=1$. If $b = 0$ or $f = 0$, then both $e$ and $c$ would have to be non-zero to ensure the determinant one condition. In this case, we can multiply a matrix on the right to add some multiple of column $2$ to column $3$ to obtain a matrix $E$ with its middle-right entry being non-zero.

$$ E = \begin{pmatrix}
0&b&c\\
0&e&f\\
1&0&0\\
\end{pmatrix}
\begin{pmatrix}
1&0&0\\
0&1&l\\
0&0&1\\
\end{pmatrix} = 
\begin{pmatrix}
0&b&c+bl\\
0&e&f+el\\
1&0&0\\
\end{pmatrix}$$\\

We can re-write $E$ in the form \\

$$\begin{pmatrix}
0&b&e\\
0&c&f\\
1&0&0\\
\end{pmatrix}$$\\

where $f \neq 0$ and $\mathrm{det}(E) = bf - ce = 1$. If $b \neq 0$, then let $F = E$. Otherwise, if $b = 0$ then we again have $c$ and $e$ are non-zero. Here, we can multiply a matrix on the left to add some multiple of row $2$ to row $1$ to obtain a matrix $F$ with its top-middle entry being non-zero.

$$ F = \begin{pmatrix}
1&k&0\\
0&1&0\\
0&0&1\\
\end{pmatrix}
\begin{pmatrix}
0&b&c\\
0&e&f\\
1&0&0\\
\end{pmatrix} = 
\begin{pmatrix}
0&b+ek&c+fk\\
0&e&f\\
1&0&0\\
\end{pmatrix}$$\\

We can re-write $F$ in the form \\

$$\begin{pmatrix}
0&b&e\\
0&c&f\\
1&0&0\\
\end{pmatrix}$$\\

where both $b$ and $f$ are non-zero. Then, the following right-multiplication can give us the matrix $G$.
$$
G = 
\begin{pmatrix}
0&b&c\\
0&e&f\\
1&0&0\\
\end{pmatrix}
\begin{pmatrix}
1&0&0\\
0&1&\frac{-c}{b}\\
0&\frac{-e}{f}&1\\
\end{pmatrix}=
\begin{pmatrix}
0&b-\frac{ce}{f}&0\\
0&0&f - \frac{ec}{b}\\
1&0&0\\
\end{pmatrix}$$\\

As the determinant one condition is preserved, then $b-\frac{ce}{f}$ and 
$f - \frac{ec}{b}$ are non-zero. We can re-write $G$ in the form \\

$$
G = 
\begin{pmatrix}
0&b&0\\
0&0&f\\
1&0&0\\
\end{pmatrix} = 
\begin{pmatrix}
0&b&0\\
0&0&\frac{1}{b}\\
1&0&0\\
\end{pmatrix}
$$\\

where above we use the fact that $b$ and $f$ are non-zero and $\mathrm{det}(G) = bf = 1$. A final right-multiplication below yields the generator $Y$.

$$\begin{pmatrix}
0&b&0\\
0&0&\frac{1}{b}\\
1&0&0\\
\end{pmatrix}
\begin{pmatrix}
1&0&0\\
0&\frac{1}{b}&0\\
0&0&b\\
\end{pmatrix}=
\begin{pmatrix}
0&1&0\\
0&0&1\\
1&0&0\\
\end{pmatrix}
= Y
$$\\

\subsubsection{Generators Proof for $Z$}

As in the previous proof, write $A$ as 

$$A=\begin{pmatrix}
a&b&c\\
d&e&f\\
g&h&i\\
\end{pmatrix}$$
\\

Again, note that at least one of $d$ and $g$ are non-zero. this allows us to left-multiply by a matrix to yield matrix $B$ with its top-left entry being zero, as shown below.

$$
B = 
\begin{pmatrix}
1&j&k\\
0&1&0\\
0&0&1\\
\end{pmatrix}
\begin{pmatrix}
a&b&c\\
d&e&f\\
g&h&i\\
\end{pmatrix}=
\begin{pmatrix}
0&l&m\\
d&e&f\\
g&h&i\\
\end{pmatrix}$$\\

If $d \neq 0$, then let $C = B$. Otherwise, then $g$ is non-zero, so left-multiply by a matrix which adds a multiple of row $3$ to row $2$ to yield matrix $C$ with its left-middle entry non-zero.

$$
C = 
\begin{pmatrix}
1&0&0\\
0&1&j\\
0&0&1\\
\end{pmatrix}
\begin{pmatrix}
a&b&c\\
d&e&f\\
g&h&i\\
\end{pmatrix}=
\begin{pmatrix}
0&l&m\\
d+jg&e+jh&f+ji\\
g&h&i\\
\end{pmatrix}$$\\

We can re-write $C$ in the form \\

$$\begin{pmatrix}
0&b&c\\
d&e&f\\
g&h&i\\
\end{pmatrix}$$\\
 
where $d$ is non-zero. Now, right-multiply by a matrix that turns the second and third entries in the middle column to zero using the first column. Let the new matrix be $D$.

$$
D = 
\begin{pmatrix}
0&b&c\\
d&e&f\\
g&h&i\\
\end{pmatrix}
\begin{pmatrix}
1&j&k\\
0&1&0\\
0&0&1\\
\end{pmatrix}=
\begin{pmatrix}
0&b&c\\
d&0&0\\
g&n&p\\
\end{pmatrix}$$\\

Since $\mathrm{det}(D) = bdp - cdn = 1$, then either $b$ or $c$ must be non-zero. If $b \neq 0$, then let $E = D$. Otherwise, we have $b = 0$ and $c \neq 0$. In this case, we can right-multiply by a matrix that adds a multiple of column $3$ to column $2$ to obtain a matrix $E$ with its top-middle entry non-zero, as shown below. \\

$$
E = 
\begin{pmatrix}
0&b&c\\
d&0&0\\
g&n&p\\
\end{pmatrix}
\begin{pmatrix}
1&0&0\\
0&1&0\\
0&j&1\\
\end{pmatrix}=
\begin{pmatrix}
0&b+cj&c\\
d&0&0\\
g&n+pj&p\\
\end{pmatrix}$$\\

We can re-write $E$ in the form \\

$$
E = 
\begin{pmatrix}
0&b&c\\
d&0&0\\
g&h&i\\
\end{pmatrix}$$\\

where $b$ is non-zero. Now we can right-multiply by a matrix that adds a multiple of column $2$ to column $3$ to obtain the matrix $F$ with its top-right entry being zero.

$$
F = 
\begin{pmatrix}
0&b&c\\
d&0&0\\
g&h&i\\
\end{pmatrix}
\begin{pmatrix}
1&0&0\\
0&1&\frac{-c}{b}\\
0&0&1\\
\end{pmatrix}=
\begin{pmatrix}
0&b&0\\
d&0&0\\
g&h&i - \frac{hc}{b}\\
\end{pmatrix}$$\\

We can re-write $F$ in the form \\

$$
F = 
\begin{pmatrix}
0&b&0\\
d&0&0\\
g&h&i\\
\end{pmatrix}$$\\

We have $\mathrm{det}(F) = -bdi = 1$, so $b$, $d$ and $i$ are non-zero. Now, right-multiply by a matrix that adds a multiple of column $3$ to columns $2$ and left-multiply by a matrix that adds a multiple of row $2$ to row $3$ to obtain the matrix $G$ with the three entries in its bottom row all equal to $i$.

$$
G = 
\begin{pmatrix}
1&0&0\\
0&1&0\\
0&j&1\\
\end{pmatrix}
\begin{pmatrix}
0&b&0\\
d&0&0\\
g&h&i\\
\end{pmatrix}
\begin{pmatrix}
1&0&0\\
0&1&0\\
0&k&1\\
\end{pmatrix}
=
\begin{pmatrix}
0&b&0\\
d&0&0\\
i&i&i\\
\end{pmatrix}
$$\\

Now left-multiply by a matrix to obtain matrix $J$ with its bottom row entries all $-1$ and left-middle entry equal to $1$.

$$
J = 
\begin{pmatrix}
-id&0&0\\
0&\frac{1}{d}&0\\
0&0&\frac{-1}{i}\\
\end{pmatrix}
\begin{pmatrix}
0&b&0\\
d&0&0\\
i&i&i\\
\end{pmatrix}=
\begin{pmatrix}
0&-bid&0\\
1&0&0\\
-1&-1&-1\\
\end{pmatrix}$$\\

Now $\mathrm{det}(J) = -(-bid)(1)(-1) = -bid = 1$, so we have 

$$
J = 
\begin{pmatrix}
0&1&0\\
1&0&0\\
-1&-1&-1\\
\end{pmatrix}
= Z
$$\\

Hence, the generator $Z$ can be attained. This completes the proof that $H$ is a maximal subgroup of $\SL$.

\end{proof} 

Indeed, consider the group $\SL$ acting on the trace 0 matrices with no eigenvectors by conjugation. The orbit of any given trace 0 matrix is its conjugacy class. By the orbit-stabiliser theorem, the size of this conjugacy class multiplied by the size of the stabiliser of this matrix is the size of the group, $\SL$, which we have already determined. Thus, counting the number of matrices in the stabiliser (which happens to be the centraliser as $X^{-1}MX = M \iff MX = XM$) of a chosen conjugacy class representative is equivalent to counting the number of elements in the conjugacy class. 

\subsection{Number of solutions to determinant 1 and trace 1 matrices stabiliser }
The number of solutions to the following equation in modulo $7$ must be computed. 
   
    $a^3+b^3+d^3 - 3a^2 d - 3ab^2 + 2db^2 - d^2 b - 3bad = 1$ \\

The equation can be written as a cubic in the variable $a$. \\

    $a^3 - 3da^2 + (-3b^2 - 3bd)a + (b^3 + d^3 + 2db^2 - d^2 b - 1)= 0$ \\

The table below shows the polynomial in $a$ for all $49$ possible pairs of $b, d$ and the corresponding possible values of $a$.
\quad
\begin{center}
    \begin{tabular}{|c|c|c|c|}
    \hline
        \textbf{Value of $b$} & \textbf{Value of $d$} & \textbf{Cubic in $a$} & \textbf{Solutions to $a$} \\
    \hline
         0 & 0 & $a^3+6$ & 1, 2, 4 \\
    \hline
         0 & 1 & $a^3+4a^2$ & 0, 3 \\
    \hline
         0 & 2 & $a^3+a^2$ & 0, 6 \\
    \hline
         0 & 3 & $a^3+5a^2+5$ & 3 \\
    \hline
         0 & 4 & $a^3+2a^2$ & 0, 5 \\
    \hline
         0 & 5 & $a^3+6a^2+5$ & 5 \\
    \hline
         0 & 6 & $a^3+3a^2+5$ & 6 \\
    \hline
         1 & 0 & $a^3+4a$ & 0 \\
    \hline
         1 & 1 & $a^3+4a^2+a+2$ & 2 \\
    \hline
         1 & 2 & $a^3+a^2+5a+1$ & \\
    \hline
         1 & 3 & $a^3+5a^2+2a+3$ & 2 \\
    \hline
         1 & 4 & $a^3+2a^2+6$ & 0, 2, 3 \\
    \hline
         1 & 5 & $a^3+6a^2+3a+5$ & 6 \\
    \hline
         1 & 6 & $a^3+3a^2+3$ & 1, 5 \\
     \hline
         2 & 0 & $a^3+2a$ & 0 \\
    \hline
         2 & 1 & $a^3+4a^2+3a$ & 0, 4, 6 \\
    \hline
         2 & 2 & $a^3+a^2+4a+2$ & 4 \\
    \hline
         2 & 3 & $a^3+5a^2+5a+5$ & 5 \\
    \hline
         2 & 4 & $a^3+2a^2+6a+1$ &  \\
    \hline
         2 & 5 & $a^3+6a^2+3$ & 2, 3 \\
    \hline
         2 & 6 & $a^3 + 3a^2+a+3$ & 4 \\
     \hline
         3 & 0 & $a^3+a+5$ & 1, 3 \\
    \hline
         3 & 1 & $a^3+4a^2+6a$ & 0 \\
    \hline
         3 & 2 & $a^3+a^2+4a+2$ & 4 \\
    \hline
         3 & 3 & $a^3+5a^2+2a+3$ & 2 \\
    \hline
         3 & 4 & $a^3+2a^2+2$ & 4 \\
    \hline
         3 & 5 & $a^3+6a^2+5a+5$ &  \\
    \hline
         3 & 6 & $a^3+3a^2+3a+4$ & \\
     \hline
         4 & 0 & $a^3+a$ & 0 \\
    \hline
         4 & 1 & $a^3+4a^2+3a+1$ &  \\
    \hline
         4 & 2 & $a^3+a^2+5a$ & 0, 1, 5 \\
    \hline
         4 & 3 & $a^3+5a^2+3$ & 4, 6 \\
    \hline
         4 & 4 & $a^3+2a^2+2a+2$ & 1 \\
    \hline
         4 & 5 & $a^3+6a^2+4a+3$ & 1 \\
    \hline
         4 & 6 & $a^3+3a^2+6a+5$ & 3 \\
     \hline
         5 & 0 & $a^3+2a+5$ & 4, 5 \\
    \hline
         5 & 1 & $a^3+4a^2+a+2$ & 2 \\
    \hline
         5 & 2 & $a^3+a^2+2$ & 2 \\
    \hline
         5 & 3 & $a^3+5a^2+6a+4$ &  \\
    \hline
         5 & 4 & $a^3+2a^2+5a$ & 0 \\
    \hline
         5 & 5 & $a^3+6a^2+4a+3$ & 1 \\
    \hline
         5 & 6 & $a^3+3a^2+3a+5$ & \\
    \hline
         6 & 0 & $a^3+4a+5$ & 2, 6 \\
    \hline
         6 & 1 & $a^3+4a^2+2$ & 1 \\
    \hline
         6 & 2 & $a^3+a^2+3a$ & 0 \\
    \hline
         6 & 3 & $a^3+5a^2+6a+5$ & \\
    \hline
         6 & 4 & $a^3+2a^2+2a+2$ & 1 \\
    \hline
         6 & 5 & $a^3+6a^2+5a+4$ & \\
    \hline
         6 & 6 & $a^3+3a^2+a+3$ & 4 \\
    \hline
    \end{tabular}
\end{center}

Altogether, we get $57$ possible triplets $(a, b, c)$ and this is the number of matrices that commute with the given matrix which had trace $0$ and determinant $1$.


\section{Background notes}
We put here more background knowledge of a general nature, much of which would be encountered in undergraduate or graduate courses. We still try to keep this parallel to the main article.

\subsubsection{Groups}

The action of group $G$ on $X$ is transitive if $X$ is non-empty and for each pair $x, y$ in $X$ there exists a $g$ in $G$ such that $g \cdot x = y$. \\

Consider a group $G$ acting on a set $X$. The orbit of an element $x$ in $X$ is the set of elements in $X$ to which $x$ can be moved by the elements of $G$. The orbit of $x$ is denoted by $G \cdot x$: 

\[ G \cdot x = \{ g \cdot x: g \in G \}. \]

The set of orbits in $X$ under action of $G$ form a partition of $X$. \\

Given $g$ in $G$ and $x$ in $X$ with $g \cdot x = x$, it is said that "$x$ is a fixed point of $g$" or that "$g$ fixes $x$". For every $x$ in $X$, the stabiliser subgroup of $G$ with respect to $x$ is the set of all elements in $G$ that fix $x$: 

\[ G_x = \{ g \in G : g \cdot x = x \}. \]

This is a subgroup of $G$ (not necessarily normal). \\

There is a bijection between the set $G / G_x$  of cosets for the stabiliser subgroup and the orbit $G \cdot x$, which is the orbit stabiliser theorem. If $G$ is finite then we also have 

\[ | G \cdot x | = |G| / |G_x| \]

\subsubsection{Matrices}

A square matrix $A$ is called diagonalisable if there exists an invertible matrix $P$ and a diagonal matrix $D$ such that 
$P^{-1}AP = D$, or equivalently $A = PDP^{-1}$. The column vectors of $P$ form a basis consisting of eigenvectors of $A$ and the diagonal entries of $D$ are the corresponding eigenvalues ot $T$ with respect to this eigenvector basis. One can raise a diagonal matrix $D$ to a power by simply raising the diagonal entries to that power and the determinant of a diagonal matrix is simply the product of all diagonal entries, such computations generalise easily to $A = PDP^{-1}$. \\


An $n \times n$ matrix $A$ over a field $F$ is diagonalisable if and only if there exists a basis of $F$ consisting of eigenvectors of $A$. If such a basis has been found, one can form the matrix $P$ having these basis vectors as columns, and $P^{-1}AP$ will be a diagonal matrix whose diagonal entries are the eigenvalues of $A$. \\

An $n \times n$ matrix $A$ is diagonalisable over the field $F$ if it has $n$ distinct eigenvalues in $F$ (i.e. if its characteristic polynomial has $n$ distinct roots in $F$), but the converse is false. E.g. if the eigenvalues are 
$1, 2, 2$ (double root of $2$) then it could be possible that the eigenspace of $A$ associated with eigenvalue $2$ has dimension $2$ (solutions to $Av=\lambda v$ when $\lambda$ is set to $2$). \\

Let $A$ be a matrix over $F$. If it is diagonalisable, then so is any power of it because $A = PDP^{-1}$ implies that 
$A^n = PD^{n}P^{-1}$ where $D^n$ is still a diagonal matrix. \\

\subsubsection{Conjugacy Classes}
In group theory, two elements $a, b$ of a group are conjugate if there is an element $g$ in the group such that $b=g^{-1}ag$. This is an equivalence relation whose equivalence classes are called conjugacy classes. \\

Members of the same conjugacy class cannot be distinguished using only group structure and therefor share many properties. \\

The \textbf{class number} of $G$ is the number of distinct conjugacy classes. All elements belonging to the same conjugacy class have the \textbf{same order}. \\

\begin{itemize}
    \item The identity element is always the only element in its class.
    \item If two elements $a, b$ belong to the same conjugacy class then they have the same order. More generally, every statement about $a$ can be translated into a statement about $b=gag^{-1}$ because the map $x \rightarrow gxg^{-1}$ is an automorphism of $G$.
    \item If $a$ and $b$ are conjugate, so are their powers $a^k$ and $b^k$. Thus, taking $k$th powers gives a map on conjugacy classes and one may conside which conjugacy classes are in its preimage. 
    \item Any element $a \in G$ lies in the centre $Z(G)$ of $G$ if and only if its conjugacy class has only one element, $a$ itself. Recall that the centre is the set of elements that commute with every element of $G$.
    \item More generally, if $C_G(a)$ is the centraliser of $a \in G$ (the subgroup consisting of all elements $g$ such that $ga=ag$), then the index $[G: C_G(a)]$ is equal to the number of elements in the conjugacy class of $a$ by the orbit-stabiliser theorem. Index refers to number of left / right cosets (same as size of group divided by size of centraliser).
    \begin{itemize}
        \item This is because we can make the G-Set $X$ be a conjugacy class of $G$. Then for any $g$ in the entire group, defined $g \cdot x = g^{-1}xg$ where $x$ is in $X$.
        \item The group action is transitive as $X$ is by definition an entire conjugacy class so there is a single orbit of any $x \in X$ whose size is the size of conjugacy class.
        \item The stabiliser subgroup of some fixed $x \in X$ is all $g \in G$ such that $x = g^{-1}xg$ or $gx=xg$ (so all $g$ that commute with $x$).
        \item Since orbit size equals conjugacy class size then by orbit-stabiliser theorem, size of conjugacy class equals size of group divided by size of stabiliser subgroup.
        \item This can be applied to $\PSL$ conjugacy classes.
    \end{itemize}
\end{itemize}

Standard geneartors of $L_3(7)$ are $a, b$ where $a$ has order $2$, $b$ has order $3$, $ab$ has order $19$ and $ababb$ has order $6$. \\

\subsubsection{Calculating entire group size} 
Using the orbit stabiliser theorem. Consider the vector 

$v = \begin{pmatrix}
1\\
0\\
0\\
\end{pmatrix}
$.

Let the set $X$ be all $7^3 - 1$ column non-zero vectors. Let $M$ be a matrix in the entire group. Let $w$ be an element of $X$. Note that $Mv = w$ whenever $M$ has its first column the same as $w$ (we can always find such an $M$ with determinant one).

This means that the orbit of $v$ is the entire set $X$. 

The stabiliser of $v$ is all matrices $M$ of the form

$$M = \begin{pmatrix}
1&a&b\\
0&c&d\\
0&e&f\\
\end{pmatrix}$$

with determinant $cf-de = 1$. So to find stabiliser size find the number of $2 \times 2$ matrices (bottom-right) which have determinant $1$. 

Can do similar computation (inductively e.g. with orbit stabiliser theorem) to get that the number of $2 \times 2$ matrices with determinant $1$ is $7\times (7^2-1)$. 

We can freely choose $a, b$ in $M$ so now we have $7^2\times 7 \times (7^2 -1)$ ways for stabiliser size. Do orbit size times stabiliser size to get $(7^3-1)(7^3)(7^2-1)$ ways
















Q



\subsection{Remaining items copied}


Maximal subgroups 

Matrices with exactly one eigenvalue
- what are the conjugacy classes of matrices whose eigenvalues are all 1.

Matrices with two distinct eigenvalues
-which conjugacy classes of matrices have eigenvalues two -1s and one 1.

generators a order 2 b order 3


\section*{Acknowledgements} We are very grateful for the encouragement and support of the parents of our talented young students many of whom are still in high school.
We would also like to acknowledge the many helpful online resources such as wikipedia (esp groupprops subwiki), ATLAS, open access initiative, Mathoverflow, mathstackexchange, overleaf, discord, youtube

\begin{enumerate}
    \item \url{http://brauer.maths.qmul.ac.uk/Atlas/v3/lin/L37/}
    \item INVERSE GALOIS PROBLEM FOR SMALL SIMPLE GROUPS DAVID ZYWINA
    \item \url{https://bit.ly/3uhb7J1}
    \item \url{https://groupprops.subwiki.org/w/index.php?title=Element-structure-special-linear-group-of-degree-three-over-a-finite-field&mobile-action=toggleviewmobile}

\end{enumerate}

{\sc Dr Michael Sun's School of Maths, NSW, Australia}

Corresponding Author: Yasiru Jayasooriya

\email{yasiruj2000@gmail.com}


\begin{thebibliography}{10}
 \bibitem{SL3gen}
  {\sc M.~Conder, E.~Robertson, P.~Williams}
  {\em Presentations for 3-Dimensional Special Linear Groups Over Integer Rings}
  Proceedings of the American Mathematical Society {\bf 115} (1992), No.1, 19-20.
 


 \end{thebibliography}
\end{document}